\documentclass[12pt]{amsart}
\usepackage{amsfonts, amssymb, amsmath, amscd, bbold, graphics, rotating}
\usepackage[all]{xy}

\providecommand{\bysame}{\leavevmode\hbox to3em{\hrulefill}\thinspace}
\providecommand{\MR}{\relax\ifhmode\unskip\space\fi MR }

\providecommand{\href}[2]{#2}

\newcommand{\C}{\mathbb{C}}
\newcommand{\tr}[1]{\mathrm{tr}(#1)}
\newcommand{\xb}{\mathbf{X}}
\newcommand{\yb}{\mathbf{Y}}
\newcommand{\zb}{\mathbf{Z}}
\newcommand{\pol}{\mathrm{pol}}
\newcommand{\ub}{\mathbf{U}}
\newcommand{\vb}{\mathbf{V}}
\newcommand{\wb}{\mathbf{W}}

\newcommand{\aq}{/\!\!/}
\newcommand{\X}{\mathfrak{X}}
\newcommand{\Y}{\mathfrak{Y}}
\newcommand{\R}{\mathfrak{R}}

\newcommand{\SL}{\mathrm{SL}(3,\C)}

\newcommand{\F}{\mathtt{F}}
\newcommand{\G}{\mathfrak{G}}
\newcommand{\xt}{\mathtt{x}}

\newcommand{\wt}{\mathtt{w}}

\newcommand{\id}{\mathbf{I}}

\newcommand{\N}[1]{\mathrm{N}_r(#1)}
\newcommand{\SLm}[1]{\mathrm{SL}(#1,\C)}
\newcommand{\GLm}[1]{\mathrm{GL}(#1,\C)}
\newcommand{\glm}[1]{\mathfrak{gl}(#1,\C)}
\newcommand{\slm}[1]{\mathfrak{sl}(#1,\C)}

\newcommand{\qt}[1]{#1^{\times r}\aq \SL}
\newcommand{\T}{\mathrm{t}}

\newtheorem{thm}{Theorem}
\newtheorem{defn}[thm]{Definition}
\newtheorem{lem}[thm]{Lemma}
\newtheorem{cor}[thm]{Corollary}
\newtheorem{rem}[thm]{Remark}
\newtheorem{prop}[thm]{Proposition}

\title[Minimal Affine Coordinates]{Minimal Affine Coordinates for $\SL$ Character Varieties of Free Groups}
\author[S. Lawton]{Sean Lawton}
\address{Departamento de Matem\'atica, Instituto Superior T\'ecnico, Lisboa, Portugal}
\address{E-mail address: $\mathtt{slawton@math.ist.utl.pt}$}
\address{URL: $\mathtt{http://www.math.ist.utl.pt/ \sim slawton}$}
\date{\today}

\begin{document}

\bibliographystyle{amsplain}

\maketitle

\begin{abstract}
Let $\X_r$ be the moduli of $\SL$ representations of a rank $r$ free group.  In this paper we determine minimal generators of the coordinate ring of $\X_r$.  This at once gives explicit global coordinates for the moduli and determines the dimension of the moduli's minimal affine embedding.  Along the way, we utilize results concerning the moduli of $r$-tuples of matrices in $\glm{3}$.  Consequently, we also state general invariant theoretic correspondences between the coordinate rings of the moduli of $r$-tuples of elements in $\glm{n}$, $\slm{n}$, and $\SLm{n}$.
\end{abstract}


\section{Introduction}
The purpose of this paper is to {\it construct} minimal affine embeddings of $\SL$-character varieties of free groups.

Let $\F_r$ be a free group of rank $r$.  The Lie group $\SL$ acts by conjugation on the space of group homomorphisms $\R_r=\mathrm{Hom}(\F_r,\SL)$.  For any such homomorphism $\rho$ let $[\rho]:=\{g\rho g^{-1}\ | \ g\in \SL\}$ be its orbit, and define the following equivalence relation in these terms: $\rho_1\sim\rho_2$ if and only if $[\rho_1]\cap [\rho_2]\not= \emptyset$.  Call the set of equivalence classes of this relation $\X_r$, and denote such a class by $\overline{[\rho]}$.  We will see $\X_r$  is an irreducible affine algebraic set (an affine variety).  Let $N_r=\frac{r}{240}\left(396+65r^2-5r^3+19r^4+5r^5\right).$  We explicitly construct subsets $\mathcal{W}_r:=\{\wt_1,...,\wt_{N_r}\}\subset \F_r$, and show

\begin{thm}
The mappings $t_{\mathcal{W}_r}:\X_r\to \C^{N_r}$ given by $$\overline{[\rho]}\mapsto (\tr{\rho(\wt_1)},...,\tr{\rho(\wt_{N_r})})$$ are polynomial embeddings where $N_r$ is minimal among all such embeddings.
\end{thm}

Most of the work in establishing this theorem is in the construction of the sets $\mathcal{W}_r$ and in establishing minimality.  We remark that given any competing set in the group ring of $\F_r$, denoted $\mathcal{S}_r\subset \C\F_r$, so that $t_{\mathcal{S}_r}$ is likewise an affine embedding, our constructive proof implicitly gives polynomial transformations between the isomorphic images $t_{\mathcal{W}_r}(\X_r)$ and $t_{\mathcal{S}_r}(\X_r)$.

To prove this theorem, we first show that only traces of evaluations at words in $\F_r$ are necessary to construct such embeddings.  This is the content of Sections 1 and 2.  In particular, the remainder of this section is devoted to setting the terms, notation, and background necessary for our discussion.  We then state our main theorem explicitly and prove its corollaries.  Section 2 discusses general relationships between $\X_r$ and some related varieties for which there are established results we find useful.  We prove our main theorem constructively in section 3.  First we show that only a subset of certain types of words in $\F_r$ are necessary (first reductions), then we count exactly how many of each type is necessary (second reductions).  Minimality will then follow from our general considerations in Section 2.

\subsection{Quotient Varieties}
Let $\G=\SLm{n}$, and $Y$ be a $\G$-variety; that is, a variety for which $\G$ acts rationally ($Y\times \G\to Y$ is regular).  The representation variety $\R_r=\mathrm{Hom}(\F_r,\G)\approx \G^{\times r}$, and the spaces $\glm{n}^{\times r}$ and $\slm{n}^{\times r}$ are affine $\G$-varieties.  $\G$ acts on each by simultaneous conjugation in each factor.  Explicitly, this action is defined as follows. Let $(A_1,...,A_r)$ be in one of $\G^{\times r}, \glm{n}^{\times r}, \text{ or } \slm{n}^{\times r}$ and let $g$ be in $\G$.  Then $$g\cdot (A_1,...,A_r)=(gA_1g^{-1},...,gA_rg^{-1}).$$ The orbit space $Y/\G$ is not generally a variety; not Hausdorff either.  However, there is a categorical quotient $\X=Y\aq \G$.  This quotient is constructed as follows.  Let $\C[Y]$ be the coordinate ring of $Y$; that is, the ring of polynomial functions on $Y$.  The conjugation action extends to $\C[Y]$ and the subring of invariants of this action, $\C[Y]^\G$, is the set of polynomial functions on $Y$ invariant under conjugation.  In other words, these polynomials are defined on orbits. But they do not distinguish orbits whose closures intersect (polynomials are continuous!).

The following definition makes clear a technical condition we will need.

\begin{defn}
A linear algebraic group $\G$ is called {\it linearly reductive} if for any rational representation $\rho:\G\to \mathrm{GL}(V)$ and any invariant vector $v\not=0$ there exists a linear invariant 
function $f$ so $f(v)\not=0$.
\end{defn}

The ``unitary trick'' shows that $\GLm{n}$ and $\SLm{n}$ are linearly reductive.

Since $\G$ is (linearly) reductive, $\C[Y]^\G$ is a finitely generated domain and so $\X=\mathrm{Spec}_{max}(\C[Y]^\G)$ is an affine variety (see Nagata \cite{Na}).

When $Y=\R_r$ the quotient $\X_r$ is called the $\G$-\emph{character variety of} $\F_r$ since it is the largest variety which parametrizes conjugacy classes of representations (characters).  In this case, there is a one-to-one correspondence between the points of $\X_r$ and the orbits of {\it completely reducible} representations (representations that are sums of irreducible representations); these are the points whose orbits are closed.  Any representation can be continuously and conjugate-invariantly deformed to one that is completely reducible, so the points of $\X_r$ are unions of orbits of representations that are deformable in this way.  Such a union is called an {\it extended orbit equivalence class}.  The character variety $\X_r$ may be accurately thought of as either the usual orbit space of $\R_r$ with the non-completely reducible representations removed, or as the the usual orbit space with extended orbit equivalences.  Either way, the resulting space is, or is in one-to-one correspondence with, an affine algebraic set, irreducible and singular, that satisfies the diagrammatic requirements needed to be a categorical quotient.  We state this definition for completeness.
\begin{defn}
Let $\G$ be an algebraic group.  A {\it categorical quotient} of a $\G$-variety $Y$ is an object $Y\aq \G$ and a $\G$-invariant morphism $\pi:Y\to Y\aq \G$ such that the following commutative diagram exists 
uniquely for all invariant morphisms $f:Y\to Z$: $$ \xymatrix{
Y \ar[rr]^-{\pi} \ar[dr]_{f} & & Y\aq \G \ar@{.>}[dl]\\
& Z &} $$
It is a {\it good} categorical quotient if the following conditions additionally hold:
\begin{enumerate}
\item[$(i)$] for open subsets $U\subset Y\aq \G$, $\C[U]\approx \C[\pi^{-1}(U)]^{\G}$;
\item[$(ii)$] $\pi$ maps closed invariant sets to closed sets;
\item[$(iii)$] $\pi$ separates closed invariant sets.
\end{enumerate}
\end{defn}
In \cite{Do} it is shown that $Y\to\mathrm{Spec}_{max}(\C[Y]^\G)$ is a good categorical quotient, if $\G$ is reductive, and so all such quotients considered in this paper are {\it good}.

\subsection{The $i^{\text{th}}$ fundamental theorem}

Perhaps our strongest motivation and our foundation to work from is the work of C. Procesi from $1976$.

In \cite{P1} Procesi shows

\begin{thm}[$1^{\text{st}}$ Fundamental Theorem of $n\times n$ Matrix Invariants]\hfill\label{prfd1}
Any polynomial invariant of $r$ matrices $A_1,...,A_r$ of size $n\times n$ is a polynomial in the invariants $\tr{A_{i_1}A_{i_2}\cdots A_{i_j}}$;  where $A_{i_1}A_{i_2}\cdots A_{i_j}$ run over
all possible noncommutative monomials.
\end{thm}

This theorem can be recast in the language used above:  $\C[\glm{n}^{\times r}\aq \G]$ is generated by traces of words in {\it generic matrices}.  Let $\C[x^k_{ij}]$ be the complex polynomial ring in $rn^2$ variables ($1\leq k \leq r$ and $1\leq i,j \leq n$), then $$\xb_k=\left(
\begin{array}{cccc}
x^k_{11} & x^k_{12} & \cdots & x^k_{1n}\\
x^k_{21} & x^k_{22}  &\cdots & x^k_{2n}\\
\vdots &\vdots & \ddots & \vdots\\
x^k_{n1} & x^k_{n2} &\cdots &x^k_{nn}\\
\end{array}\right)$$ are generic matrices.  For any word $\wt$ in $\F_r=\left< \xt_1,...,\xt_r\right>$, let $\wb$ be the word $\wt$ with each $\xt_i$ replaced by the generic matrix $\xb_i$.  In these 
terms, the first fundamental theorem says $$\C[\glm{n}^{\times r}\aq \G]=\C[\tr{\wb}\ |\ \wt\in \F_r].$$   

Procesi showed in \cite{P1} that the index $j$ in Theorem $\ref{prfd1}$ is bounded:  $j\leq 2^n-1$.  It is called the {\it degree of nilpotency} and is often denoted $d(n)$.  In $1974$ Razmyslov \cite{Ra} had shown that $j\leq n^2$.  For $1\leq n\leq 4$, it is known that $j=\frac{n(n+1)}{2}$ (and conjectured to be true in general).  See \cite{DF} for more on the Kuzmin conjecture.  Consequently, the word length $|\wt|$ is bounded; not at all obvious.

A first fundamental theorem in invariant theory describes sufficient generators, a second fundamental theorem describes sufficient relations.

\begin{thm}[$2^{\text{\tiny nd}}$ Fundamental Theorem of $n\times n$ Matrix Invariants]\hfill \label{prfd2}
Let $\chi(t)=\sum c_k(\xb)t^{n-k}$ be the formal expression of the characteristic polynomial $\mathrm{det}(\xb-t\id)$.  Then any relation in $\C[\glm{n}^{\times r}\aq \G]$ is a consequence of the formal expression of $\tr{\chi(\xb)\cdot\xb}$.
\end{thm}

From Theorem $\ref{prfd1}$ any relation in $\C[\glm{n}^{\times r}\aq \G]$ is necessarily a polynomial expression whose terms are products of the invariants $\{\tr{A_{i_1}A_{i_2}\cdots A_{i_j}}\}$.  In Theorem $\ref{prfd2}$ the word ``consequence'' is vague, but can be made precise.  A polynomial identity $g=0$ is a {\it consequence of the polynomial identities} $f_i=0,\ i\in I$ if any algebra satisfying the identities $f_i=0$ also satisfies $g=0$.  The consequences of the expression $\tr{\chi(\xb)\cdot\xb}$ come from its multi-linearizations.  See \cite{DF} and \cite{P1} for further details.

Quoting Procesi:

\begin{quote}According to the general theory, we will split the description into two steps.  The so called ``first fundamental theorem,'' i.e., a list of
generators for $T_{i,n}$, and the ``second fundamental theorem,'' i.e., a list of relations among the previously found generators.  Of course, it would be very
interesting to continue the process by giving the ``$i^{\text{th}}$ fundamental theorem,'' i.e., the full theory of syzigies; unfortunately, this seems to be still out of
the scope of the theory as presented in this paper.\end{quote}

This remains true enough 30 years later.

\subsection{Some Progress}

In 2003 Drensky \cite{Dr} gave a complete and uniform description of the invariant ring of $2\times 2$ matrices.  In \cite{ATZ} and \cite{T3}, algebraically independent generators are worked out for $2\times 2$ matrices.  Between $1958$-$1971$, work was done establishing minimal generators for the invariants of products of arbitrary $3\times 3$ matrices 
\cite{M,MS1,MS2,MS3,Sp1,SR1,SR2,SR3}.  However, in 1989 \cite{AP} gave an algorithm to establish minimality in general and implemented it for $3\times 3$ matrices.
In 2002 Nakamoto (see \cite{Na1}) describes the $\mathbb{Z}$-scheme of two $3\times 3$ generic matrices; and later (working over a field of characteristic $0$) \cite{ADS} also describes the ideal.  Recently, exciting new results using methods similar to those in \cite{AP} were established in \cite{BD1} concerning the ideal of relations for generic $3\times 3$ matrices.  In particular, the minimal degree of generators of the ideal of relations was found to be $7$ and the degree $7$ relations were then described in general.  We note that this is an incomplete description of contributions.  See \cite{DF} for a more thorough account.

For two {\it unimodular} $3\times 3$ matrices, in \cite{L2} we prove

\begin{thm} \label{Lawton1} Let $\X=\SL^{\times 2} \aq \SL$.  Then the following hold:
\begin{enumerate}
\item[$(i)$] $\C[\X]$ is minimally generated by the nine affine coordinate functions
\begin{align*}
\mathcal{G}=&\{\tr{\xb_1},\tr{\xb_2},\tr{\xb_1\xb_2},\tr{\xb_1^{-1}},\tr{\xb_2^{-1}},\tr{\xb_1\xb_2^{-1}},\\ &\tr{\xb_2\xb_1^{-1}},\tr{\xb_1^{-1}\xb_2^{-1}},\tr{\xb_1\xb_2\xb_1^{-1}\xb_2^{-1}}\}.
\end{align*}
\item[$(ii)$] The eight elements in $\mathcal{G}\backslash \{\tr{\xb_1\xb_2\xb_1^{-1}\xb_2^{-1}}\}$ are a maximal algebraically independent subset.  Therefore, they are local parameters, since the Krull dimension of $\X$ is $8$.
\item[$(iii)$] $\tr{\xb_1\xb_2\xb_1^{-1}\xb_2^{-1}}$ satisfies a monic (degree 2) relation over the algebraically independent generators.  It generates the ideal.
\item[$(iv)$] $\mathrm{Out}(\F_2)$ acts on $\C[\X]$ and has an order $8$ subgroup which acts as a permutation group on the independent generators; as such distinguishes them.
\end{enumerate}
\end{thm}

This paper marks our first step to generalize this theorem to the general case of $\SL^{\times r} \aq \SL$ for arbitrary values of $r$.

\begin{rem}
In Chapter 10 of \cite{Fo}, parts $(i)$ and $(iii)$ of Theorem \ref{Lawton1} were earlier established.  Although an explicit formula is not derived by Fogg, he provides an exact means to compute the ideal.  Compare these results also to results in \cite{We,Si2}. 
\end{rem}

\subsection{Main Results}

Our main theorem is 

\begin{thm}\label{maintheorem}
$\C[\SL^{\times r}\aq \SL]$ is minimally generated by 
$\binom{r}{1}$ invariants of the form $\tr{\xb}$, 
$\binom{r}{1}$ invariants of the form $\tr{\xb^{-1}}$, 
$\binom{r}{2}$ invariants of the form $\tr{\xb\yb}$, 
$2\binom{r}{2}$ invariants of the form $\tr{\xb\yb^{-1}}$,
$\binom{r}{2}$ invariants of the form $\tr{\xb^{-1}\yb^{-1}}$,
$\binom{r}{2}$ invariants of the form $\tr{\xb\yb\xb^{-1}\yb^{-1}}$,
$2\binom{r}{3}$ invariants of the form $\tr{\xb\yb\zb}$, 
$6\binom{r}{3}$ invariants of the form $\tr{\xb\yb\zb^{-1}}$, 
$3\binom{r}{3}$ invariants of the form $\tr{\xb\yb\zb\yb^{-1}}$, 
$6\binom{r}{3}$ invariants of the form $\tr{\xb\yb^{-1}\zb^{-1}}$,
$6\binom{r}{3}$ invariants of the form $\tr{\xb\yb\zb^{-1}\yb^{-1}}$, 
$\binom{r}{3}$ invariants of the form $\tr{\xb^{-1}\yb^{-1}\zb^{-1}}$,
$5\binom{r}{4}$ invariants of the form $\tr{\wb\xb\yb\zb}$,
$20\binom{r}{4}$ invariants of the form \\$\tr{\wb\xb\yb\zb^{-1}}$,
$18\binom{r}{4}$ invariants of the form $\tr{\wb\xb\yb^{-1}\zb^{-1}}$, 
$8\binom{r}{4}$ invariants of the form $\tr{\wb\xb\yb\zb\yb^{-1}}$,
$12\binom{r}{5}$ invariants of the form $\tr{\ub\vb\wb\xb\yb}$,
$35\binom{r}{5}$ invariants of the form $\tr{\vb\wb\xb\yb\zb^{-1}}$, and
$15\binom{r}{6}$ invariants of the form $\tr{\ub\vb\wb\xb\yb\zb}$.
\end{thm}

We give a numeric and geometric consequence of our main theorem.

\begin{cor}
The number of minimal generators for $\C[\SL^{\times r}]^{\SL}$ is $$N_r=\frac{r}{240}\left(396+65r^2-5r^3+19r^4+5r^5\right).$$
\end{cor}

\begin{proof}
Adding up the number of generators from Theorem \ref{maintheorem}, we conclude the sum
\begin{align*}
&2\binom{r}{1} + 5\binom{r}{2} + 24\binom{r}{3} + 
 51\binom{r}{4} + 47\binom{r}{5} + 15\binom{r}{6} \\ &\\
&=2 r+\frac{5}{2} r(r-1)  +4 r(r-1)(r-2) +\frac{17}{8} r(r-1) (r-2) (r-3) \\
&+\frac{47}{120} r(r-1) (r-2) (r-3) (r-4)+\frac{1}{48} r(r-1) (r-2) (r-3) (r-4) (r-5)\\
&\\
&=\frac{r}{240}  \left(5 r^5+19 r^4-5 r^3+65r^2+396\right).
\end{align*}
\end{proof}

The minimal generators are coordinate functions for the variety $\X_r$; that is, letting $\{t_1,...,t_{N_r}\}$ be polynomial indeterminates over $\C$, there exists a finitely generated ideal, $\mathfrak{I}$, so $$\X_r= \mathrm{Spec}_{max}\left(\C[t_1,...,t_{N_r}]/\mathfrak{I}\right).$$  Consequently, we have

\begin{cor}\label{embedding}
There exists an affine embedding $\X_r\to\C^{N_r}$ where $N_r$ is minimal among all affine embeddings $\X_r\to \C^{N}$.
\end{cor}
 
Moreover, let $\{\wt_1,...,\wt_{N_r}\}$ be any set of $N_r$ words in $\F_r$ corresponding to a minimal generating set for $\C[\X_r]$.  Then the embedding from Corollary $\ref{embedding}$ is given by
$$\overline{[\rho]}\mapsto (\tr{\rho(\wt_1)},...,\tr{\rho(\wt_{N_r})}),$$ where $\rho$ is in $\R_r$ and $\overline{[\rho]}$ is an extended equivalence class in $\X_r$ from the projection $\R_r\to \X_r$.

There are different choices for $\{\wt_1,...,\wt_{N_r}\}$ beyond simple cyclic permutations of the letters in the words.  We will see in the derivations of the minimal generating sets (the relations for which we call {\it reduction relations}) explicit formulas to algebraically change variables and global coordinates.  Abstractly, we know that for two affine embeddings $\X_r\to \C^{N_r}$ there is a polynomial mapping $\C^{N_r}\to \C^{N_r}$ which commutes with the two embeddings.  Our explicit reduction formulas give concrete form to these polynomial mappings.  In other words, the reduction formulas {\it constructively} give examples of equivalent global coordinate systems on $\X_r$.

\begin{rem}
Having an explicit global description of $\X_r$, we want an equally explicit description of the local coordinates of an affine patch in $\X_r$.  We address this and some very interesting symmetry in $\X_r$ (group actions that preserve patches) in future work.
\end{rem}

\begin{rem}
In 2005 Lopatin \cite{Lo2}, \cite{Lo3} constructed a minimal generating set for $K[\mathfrak{gl}(3,K)^{\times r}]^{\mathrm{GL}(3,K)}$ for any infinite field $K$ of arbitrary characteristic. Using our Proposition $\ref{slgl}$ with the results in \cite{Lo3} one can likewise construct a minimal system of generators for $\C[\SL^{\times r} \aq \SL]$.
\end{rem}

\section{$\SLm{n}^{\times r}\aq \SLm{n},\  \glm{n}^{\times r}\aq\GLm{n}$, and $\slm{n}^{\times r}\aq \SLm{n}$ }\label{sec2}

From the first and second fundamental theorems of Procesi, we have coordinates for $\glm{n}^{\times r}\aq \SLm{n}$.  We now show how these coordinates relate to coordinates for the space $\SLm{n}^{\times r}\aq \SLm{n}$, which is our principal interest.

\subsection{Preserving Trace Generators}

The difference between the moduli of {\it arbitrary} $n\times n$ matrices and matrices with determinant $1$ (unimodular) is the inclusion or exclusion of the invariants of the form $\tr{\xb^n}$ in their coordinate rings.  In otherwords, we have

\begin{prop}\label{slgl}
$$\C[\SLm{n}^{\times r}\aq \SLm{n}]\approx \C[\glm{n}^{\times r}\aq \GLm{n}]/\mathfrak{I},$$ where $\mathfrak{I}=(\tr{\xb_1^n}-P(\xb_1),...,\tr{\xb_r^n}-P(\xb_r))$.

In general, $P(\xb)=$ $$\tr{\xb^n}+(-1)^nn(\mathrm{det}(\xb)-1)=(-1)^{n+1}n+\sum_{k=1}^{n-1}(-1)^{n+k+1}C_0^{n-k}(\xb)\tr{\xb^k}$$ is a polynomial in terms of the coefficients of the
characteristic polynomial.  

These coefficients may be computed recursively using $$C^0_0(\xb)=1,\   C^1_0(\xb)=\tr{\xb}, \text{and } C_0^n(\xb)=\frac{1}{n}\sum_{k=1}^n (-1)^{k+1}C_0^{n-k}(\xb)\tr{\xb^k}.$$
\end{prop}

\begin{rem}
In the case of $\SL$, $P(\xb)=3+\frac{3\tr{\xb}\tr{\xb^2}-\tr{\xb}^3}{2}$.
\end{rem}  

The proof of Proposition $\ref{slgl}$ will follow from two lemmas:  the ideal cuts out $\SLm{n}^{\times r}$ from $\glm{n}^{\times r}$ as schemes; and the ideal passes through the quotient.  After the proof, we will show how each relates to $\slm{n}^{\times r}\aq \SLm{n}$.

\begin{rem}
It is worth mentioning that although the entries $x_{ij}^k$ of the generic matrices $\xb_k$ are elements of the polynomial ring $\C[x_{ij}^k]$ with $rn^2$ indeterminates, the entries of the {\it unimodular} generic matrices are elements in $$\C[x_{ij}^k]/(\det(\xb_k)-1).$$  In particular, letting $\overline{x_{ij}^k}$ represent the projection of $x_{ij}^k$, the 
unimodular generic matrices take the form $\overline{\xb_k}=\big(\overline{x_{ij}^k}\big)$.  However, we will not distinguish in notation between these two types of generic matrices since the context is always clear.  In other words, the ``overline'' will from this point be omitted.
\end{rem}

We now proceed with the proof of Proposition \ref{slgl}.

\begin{proof}

First, we address a technical point:  the quotient is the same if we act by $\GLm{n}$ instead of $\SLm{n}$.  Indeed, for any $\xb$ in $\GLm{n}$, $\yb=\det(\xb)^{-\frac{1}{3}}\xb$ is in $\SLm{n}$.  Let $\overrightarrow{\xb}=(\xb_1,...,\xb_r)$ be in $\SLm{n}^{\times r}$, $\glm{n}^{\times r}$, or $\slm{n}^{\times r}$.  Then $$(\xb\xb_1\xb^{-1},...,\xb\xb_r\xb^{-1})=(\yb\xb_1\yb^{-1},...,\yb\xb_r\yb^{-1}).$$  In other words, the actions are identical.  Strictly speaking the matrix $\yb$ above roughly depends on an $n^{\text{th}}$ root of unity; in other words, there are up to $n$ solutions to the equation $x^n-\det(\xb)=0$.  But this poses no issue since there is always at least one solution and any such solution gives rise to an identical action.  Any way you look at it, the orbit $\GLm{n}\overrightarrow{\xb}\subset \SLm{n}\overrightarrow{\xb}$, but since $\SLm{n}\subset \GLm{n}$ we have the reverse inclusion as well; the orbits are identical.  Now since in the cases we are considering the ``extended orbits'' of the categorical quotient are determined by the orbits themselves (more precisely their closures and how they intersect), the fact that the orbit spaces are equal implies that for $Y$ equal to any of $\SLm{n}^{\times r}$, $\glm{n}^{\times r}$, or $\slm{n}^{\times r}$, we have $$Y\aq \GLm{n}=Y\aq\SLm{n}.$$

We come to the first lemma.

\begin{lem}\label{lem1}
Let $\G$ be a linearly reductive algebraic group acting rationally on a $\C$-algebra $\C[Y]$ leaving an ideal $I$ invariant $(\G I\subset I)$.  Then $$\left(\C[Y]/I\right)^\G\approx 
\C[Y]^\G/\left(I\cap \C[Y]^\G\right).$$
\end{lem}

\begin{proof}
The inclusion $\C[Y]^\G\subset \C[Y]$ induces an injection $$\C[Y]^\G/\left(I\cap \C[Y]^\G\right)\longrightarrow \left(\C[Y]/I\right)^\G.$$  The content of this 
lemma is that it is surjective.  This follows from the assumption that $\G$ is {\it linearly} reductive.  See page 43 in \cite{Do}.
\end{proof}

\begin{lem}\label{lem2}As affine $\C$-algebras,
$$\C[\R_r]\approx \C[\glm{n}^{\times r}]/(\det(\xb_1)-1,...,\det(\xb_r)-1)$$
\end{lem}

\begin{proof}
Since $\C[\R_r]\approx \C[\SLm{n}]^{\otimes r}$ it is sufficient to show that $$\C[\SLm{n}]\approx \C[\glm{n}]/(\det(\xb)-1).$$  However, $\SLm{n}$ is defined to be the solutions of $\det(\xb)-1=0$ in all $n\times n$ matrices; implying that $\R_r$ is cut-out of $\glm{n}^{\times r}$ as sets.  Moreover, the determinant is an irreducible polynomial.  Thus, $\C[\glm{n}]/(\det(\xb)-1)$ is a reduced algebra (no non-trivial nilpotents). Therefore, $\R_r$ is cut-out of $\glm{n}^{\times r}$ by the ideal $(\det(\xb_1)-1,...,\det(\xb_r)-1)$ as schemes.
\end{proof}

With these lemmas established, we return to the main argument for Proposition $\ref{slgl}$.
Let $\mathfrak{I}_R=(\det(\xb_1)-1,...,\det(\xb_r)-1)$ be the ideal generated in $R$, and let $Y=\glm{n}^{\times r}$.  Consequently, by Lemmas \ref{lem1} and \ref{lem2} $$\C[\R_r]^\G\approx 
\left(\C[Y]/\mathfrak{I}_{\C[Y]}\right)^\G\approx \C[Y]^\G/\left(\mathfrak{I}_{\C[Y]}\cap \C[Y]^\G\right)\approx \C[Y]^\G/\mathfrak{I}_{\C[Y]^\G}.$$
The last isomorphism follows since the generators of the ideal are themselves invariants and hence they are fixed by the $\G$-action (this observation also shows the ideal is stable under the 
$\G$-action; a necessary assumption).

Clearly, since we are working over a field of characteristic $0$, the ideal generated by the polynomials $\det(\xb_k)-1$ and the one generated by the polynomials 
$\tr{\xb_k^n}-P(\xb_k)=(-1)^{n+1}n(\mathrm{det}(\xb_k)-1)$ are identical.  

To finish the proof, it remains to derive the recursion formula for the polynomials $P(\xb)$.

\begin{lem} $$P(\xb)=(-1)^{n+1}n+\sum_{k=1}^{n-1}(-1)^{n+k+1}C_0^{n-k}(\xb)\tr{\xb^k},$$ where $C^0_0(\xb)=1,\ C^1_0(\xb)=\tr{\xb},$ and $$C_0^n(\xb)=\frac{1}{n}\sum_{k=1}^n 
(-1)^{k+1}C_0^{n-k}(\xb)\tr{\xb^k}.$$
\end{lem}

\begin{proof}
First let $$\det(t\id-\xb)=\sum_{k=0}^{n} (-1)^{n-k}C^n_k(\xb)t^k$$ define the coefficients of the Cayley-Hamilton polynomial for an $n\times n$ matrix $\xb$.  We know that $C^n_n(\xb)=1$, 
$C^n_{n-1}(\xb)=\tr{\xb}$ and $C^n_0(\xb)=\det(\xb)$.  By Newton's trace formulas each $C^n_k(\xb)$ is a polynomial in the traces of powers of the matrix $\xb$.  Polynomials in traces of powers of 
a matrix are functions that make sense for any size matrix, so the domain of the functions $C^n_k(\xb)$ can be extended to include any size matrix $\xb$ (not only $n\times n$, as the 
superscript originally implied).  Observing this and that $\mathrm{deg}(C^n_k(\xb))=n-k$ (the determinant is homogeneous of degree $n$), we have $$C^n_k(\xb)=C^{n-i}_{k-i}(\xb),$$ since they are also 
the elementary symmetric functions in the eigenvalues of the matrix $\xb$.  Consequently, $C^n_k(\xb)=C^{n-k}_0(\xb)$.  As already noted $C^{n-k}_0(\xb)=\det(\xb)$ for an $(n-k)\times (n-k)$ matrix.  
Solving the characteristic polynomial in this case gives the required recursion.
\end{proof}

With the above lemma complete, the proof of Proposition \ref{slgl} is likewise complete.
\end{proof}

Using Lemma \ref{lem1}, it further follows that \begin{equation}\C[\slm{n}^{\times r}\aq\G]\approx \C[\glm{n}^{\times r}\aq  \GLm{n}]/(\tr{\xb_1},...,\tr{\xb_r}).\nonumber\end{equation}  
One can say more, since by a change of generators induced by the map $$\xb\mapsto \xb-\frac{1}{n}\tr{\xb}\id$$ (see \cite{DF}), it follows that $$\C[\slm{n}^{\times 
r}\aq\G][\tr{\xb_1},...,\tr{\xb_r}]\approx\C[\glm{n}^{\times r}\aq \GLm{n}].$$  In other words, the exact sequence of $\C[Y]^{\G}$-modules $$0\to (\tr{\xb_1},...,\tr{\xb_r})\to \C[\glm{n}^{\times 
r}\aq  \GLm{n}]\to \C[\slm{n}^{\times r}\aq\G]\to 0$$ splits.  

Consequently, $$\C[\R_r\aq \G]\approx\C[\slm{n}^{\times r}\aq\G][\tr{\xb_1},...,\tr{\xb_r}]/\mathfrak{I}.$$

Using these isomorphisms, one can obtain results relative to all three quotient varieties:  $\X_r=\SLm{n}^{\times r}\aq \SLm{n}$, $\mathfrak{Z}_r=\slm{n}^{\times r}\aq \SLm{n}$, and 
$\mathfrak{Y}=\glm{n}^{\times r}\aq \GLm{n}$ from any one alone.  Our motivating interest is in $\X_r$ and so we focus our attention here.  We will however switch between $\X_r$ and $\Y_r$ as needed.

\begin{rem}
We note the Krull dimensions:  $\mathrm{dim}\X_r+r=\mathrm{dim}\Y_r=\mathrm{dim}\mathfrak{Z}_r+r$.
\end{rem}

\begin{rem}
This of course begs the question, how does $\C[\GLm{n}^{\times r}]^{\GLm{n}}$ relate to these varieties?  The answer is that $\GLm{n}^{\times r}$ is a quasi-affine variety of $\glm{n}^{\times r}$.  In 
particular, it is the principle open set defined by the product of the determinants of the generic matrices.  Since the determinant is an invariant function, and ``taking invariants'' commutes with ``localization,''  we have $$\C[\GLm{n}^{\times r}]^{\GLm{n}}\approx \C[\glm{n}^{\times r}\aq \GLm{n}]\left[\frac{1}{\det(\xb_1)\cdots\det(\xb_r)}\right], $$ where $\C[\glm{n}^{\times r}\aq \GLm{n}]\left[\frac{1}{\det(\xb_1)\cdots\det(\xb_r)}\right]$ is the localization at the product of determinants.
\end{rem}

\subsection{Preserving Minimality}

Before we continue, we prove that the number of minimal generators for $\C[\X_r]$ is exactly $r$ less than that of $\C[\Y_r]$.   We refer to the projection $$\Pi:\C[\Y_r]\longrightarrow \C[\X_r]\approx \C[\Y_r]/\mathfrak{I}$$ as ``the projection from Proposition \ref{slgl}.''  What we actually show is that $\Pi$ preserves minimality (or minimal sets push forward) when the set of generators for $\C[\Y_r]$ has the form $\{\tr{A_{i_1}A_{i_2}\cdots A_{i_j}}\}$.  We call such a set of generators {\it Procesi generators} if additionally no generator has the form $\tr{\wb_1\xb^n\wb_2}$ where at least one of the words $\wb_i$ is not the identity.  Using the characteristic polynomial $\sum c_k(\xb)\xb^{n-k}=0$ one can always arrange for any set of minimal generators of $\C[\Y_r]$ to be Procesi generators.

\begin{prop}\label{numbgen}
Let $\mathcal{G}$ be a minimal set of Procesi generators for $\C[\Y_r]$.  Then $\Pi\left(\mathcal{G}-\{\tr{\xb_1^n},...,\tr{\xb_r^n}\}\right)$ is a minimal set of generators for $\C[\X_r]$.
\end{prop}

\begin{proof}
We first show we can eliminate $r$ of the minimal generators of $\C[\Y_r]$ after projection and then we show that we cannot eliminate any more.  Any set of Procesi generators must include the set $\{\tr{\xb_1^n},...,\tr{\xb_r^n}\}$.  If not then locally the determinant $\mathrm{det}(\xb_i)$ (a polynomial only involving the matrix entries of $\xb_i$) would be expressible as a polynomial in $\tr{\xb_i^{n-k}}$ for $0<k<n$.  But since the coefficients of the characteristic polynomial are algebraically independent, there can be no such relation.

Moreover, once we assume the determinant is 1, we can freely remove the expression $\tr{\xb_i^n}$ in $\C[\X_r]$.  Said differently the structure of the ideal $\mathfrak{I}$ in Proposition $\ref{slgl}$ allows us to freely remove the $r$ generators $\tr{\xb_i^n}$ in the ring $\C[\X_r]$.

To complete the proof of Proposition $\ref{numbgen}$, it remains to prove that there are no further reductions after choosing a minimal set of generators for $\C[\Y_r]$.

If there was a further reduction after passing to unimodular invariants, then there would be a relation of the form $\tr{\wb}=Q(\xb_1,...,\xb_r)$ where $\wb$ is a word of some length corresponding to a minimal generator in $\C[\Y_r]$ and $Q$ is a polynomial trace expression in terms of generic matrices not including a term with a factor $\tr{\wb}$.  Moreover, because $\C[\X_r]$ is filtered, $Q$ can be assumed to have no term (with respect to the generic unimodular matrix entries) with degree greater than the length of $\wb$.  We additionally assume that $\wb\not=\xb^n$ since we have already eliminated these $r$ generators.

Pulling back from the projection in Proposition $\ref{slgl}$, there exists polynomial trace expressions $f_1,...,f_r\in \C[\Y_r]$ so \begin{equation}\label{eq:4}\tr{\wb}-Q(\xb_1,...,\xb_r)=\sum f_i\left(\tr{\xb_i^n}-P(\xb_i)\right)\end{equation} in the graded ring of arbitrary $n\times n$ invariants $\C[\Y_r]$.  However, the degree of the left-hand-side and the degree on the right-hand-side of Equation $\eqref{eq:4}$ must be equal, which implies $\tr{\wb}$ cannot be part of any $f_i$, unless $\tr{\wb}$ is of the form $\tr{\xb^n}$ and $f_i$'s are constant; we assumed this was not the case.  Thus we would have a further reduction in the ring of arbitrary invariants, which contradicts the minimality of the generators of $\C[\Y_r]$.
\end{proof}

Let $\N{x,y}$ be the minimal number of generators in $\C[\glm{n}^{\times r}\aq \G]$ of word length $x$ in $y$ letters.  So $1\leq x\leq d(n)$ and $y\leq x$.  As a convention, we say $\N{x,y}=0$ if $y>x$.

For example, consider $\glm{3}$.  In this case, $\mathrm{N}_3(2,1)=3$ since the only length $2$ words that may be constructed from three letters using only one letter per word are $\xb^2,\yb^2,$ and $\zb^2$.  It is reasonable to assume that in general, $\N{2,1}=r=\binom{r}{1}$.  The following combinatorial lemma proves this and a generalization of this fact as well.

\begin{prop}\label{count}
$\N{x,y}=\mathrm{N}_y(x,y)\binom{r}{y}$.
\end{prop}

\begin{proof}
Observe that since $\C[\Y_r]$ is multigraded, any reduction will arise from a homogeneous multidegree relation.  Consequently, once the number of generators in a certain degree are determined in the number of letters corresponding to that degree (words of length $l$ for $r=l$), then for $r\geq l$ the number of generator types is determined.  The multidegree is exactly determined by the total degree (word length) and the number of each type of letter in the word.  Since there cannot be relations among generators of differing multidegree, there cannot be relations among generators differing in total degree (word length) and the number of distinct letters in the word (and their multiplicities).

Before we move on let us consider another example in the case $\glm{3}$.  For $r=3$ the generators of degree $3$ in three distinct letters, multidegree $(1,1,1)$, are $\tr{\xb\yb\zb}$ and $\tr{\xb\zb\yb}$, which is immediate by cyclic equivalence.  It is not immediate that there are no further reductions, but we will assume it for now to illustrate our point (minimality results of \cite{AP} establish it).  Regardless, when $r=4$ there are four combinations of three letters and so according to this proposition there are $8$ generators of this type.  Cyclic reduction promises that there are no \emph{more} than $8$.  Suppose there was less and let $\wb$ be the fourth letter.  There are only two generators with multidegree $(1,1,1,0)$; namely, the generators $\tr{\xb\yb\zb}$ and $\tr{\xb\zb\yb}$.  But there is no relation among them from the $r=3$ case.  The same sort of reasoning can be applied to say $\tr{\xb\yb\wb}$ and $\tr{\xb\wb\yb}$; the only generators of multidegree $(1,1,0,1)$.  Again, there is no relation among them for the same reason as before, the $r=3$ case.

Returning to the argument, since the letters used in a word and the word's length completely determine its multidegree, and since there can be no relations except among generators of the same multidegree, the number of such generators is exactly the number of such generators occurring for the first time ($r=y$) times the number of possible differing multidegrees of the same type (after ignoring zeros); that is, $\binom{r}{y}$.
\end{proof}

Again for the case $\glm{3}$ we work an example which we believe motivates our next section and shows how to use these propositions.

In \cite{AP} it is shown that when $r=3$ a basis for the degree 3 generators (partitions $(3)$ and $(1,1,1)$ in the next section) is  \begin{align*}&\{\tr{\xb^3},\tr{\yb^3},\tr{\zb^3},\tr{\xb^2\yb},\tr{\xb^2\yb}, \tr{\zb^2\xb},\\ &\tr{\xb^2\zb},\tr{\yb^2\zb},\tr{\zb^2\yb}, \tr{\xb\yb\zb},\tr{\xb\zb\yb}  \}.\end{align*}  Proposition $\ref{count}$ shows for $r\geq 3$ there are $\binom{r}{1} + 2\binom{r}{2} + 2\binom{r}{3}$ generators of degree $3$; exactly $r$ of type $\tr{\xb^3}$.  However, Proposition $\ref{numbgen}$ allows us to freely remove the $r$ generators $\tr{\xb^3}$, and no others.  This leaves us with all the generators for $\C[\X_r]$ which arise from the degree 3 generators in $\C[\Y_r]$, for any value of $r$.  Continuing this process for all possible multi-degrees will provide us with a minimal generating set.

\section{$\SL$ Minimal Generators and Reductions}

\subsection{First Reduction: Generator Types}\label{firstreduc}

We begin this section by reviewing some of the results in \cite{L2}.  Let $\G=\SL$ and $\X_r=\G^{\times r}\aq \G$.  Capital bold 
letters $\ub, \vb, \wb$ will denote words in the generic matrices (unimodular) unless otherwise stated, and often $\xb, \yb, \zb$ will denote words of length 1.  We will frequently speak of word length.  When counting the length of a word, letters with a negative power are counted twice (this was called {\it weighted word length} in \cite{L2}).  For instance, $\xb^{-1}$ has length 2.  We denote the word length of a word $\wb$ by $|\wb|$.  It will be assumed, since we are ultimately concerned about traces of such words, that words have been cyclically reduced; and so length is computed after such a reduction.

We will additionally find useful the following two definitions, both in terms of (weighted) length of words in generic (unimodular) matrices.

\begin{defn}
The degree of a polynomial expression in generic matrices $($generic unimodular matrices, respectively$)$ with coefficients in $\C[\Y_r]$ $(\C[\X_r]$, respectively$)$ is the largest word length $($weighted word length, respectively$)$ of monomial words in the expression that is minimal among all such expressions.
\end{defn}

\begin{defn}
The trace degree of a polynomial expression in generic matrices $($generic unimodular matrices, respectively$)$ with coefficients in $\C[\Y_r]$ $(\C[\X_r]$, respectively$)$ is the maximal degree over all monomial words within a trace coefficient of the expression.
\end{defn}

In \cite{L2} the author shows that $\C[\X_r]$ is generated by $\{\tr{\wb}\ |\ \ |\wb|\leq 6\}$.  The Cayley-Hamilton equation provides the identity $$\xb^2-\tr{\xb}\xb 
+\tr{\xb^{-1}}\id-\xb^{-1}=0.$$  So we may freely replace any polynomial generator of the form $\tr{\ub\xb^2\vb}$ with one of the form $\tr{\ub\xb^{-1}\vb}$. This follows since $$\tr{\ub\xb^{2}\vb}=\tr{\ub\xb^{-1}\vb}+\tr{\xb}\tr{\ub\xb\vb} -\tr{\xb^{-1}}\tr{\ub \vb};$$ justifying the {\it weight} in the {\it weighted length}.  Therefore, the ring of invariants is generated by traces of words whose letters have exponent $\pm 1$, of word length 6 or less.

We just showed that for $\G=\SL$ every generator in terms of a letter with exponent $-1$ is {\it interchangeable} with a generator in terms of the exact same word with the $-1$ exponents replaced with exponents of $2$.  For instance, one can replace $\tr{\wb\xb^{-1}}$ with $\tr{\wb\xb^2}$, etc.

By linearizing the Cayley-Hamilton polynomial (see \cite{L2}), we deduce
\begin{align}\label{eq:1}
&\yb\xb^2+\xb^2\yb+\xb\yb\xb=\nonumber\\
\nonumber\\
&\tr{\xb}\big(\yb\xb +\xb\yb\big)+\big(\tr{\xb\yb}-\tr{\xb}\tr{\yb}\big)\xb+\left(\frac{\tr{\xb^2}-\tr{\xb}^2}{2}\right)\yb\nonumber\\
&+\left(\tr{\yb\xb^2}-\tr{\xb}\tr{\xb\yb}+\tr{\yb}\frac{\tr{\xb}^2 -\tr{\xb^2}}{2}\right)\id +\tr{\yb}\xb^2.
\end{align}  
Define $\pol(\xb,\yb)$ to be the right hand side of Equation \eqref{eq:1}.  Then $\pol(\xb,\yb)=\yb\xb^2 +\xb^2\yb+\xb\yb\xb.$

Multiplying on the left by a word $\wb_1$ and on the right by a word $\wb_3$, substituting a word $\wb_2$ for $\yb$, letting $\xb$ possibly be its inverse, and taking the trace of both sides of 
Equation \eqref{eq:1} yields  
\begin{align}\label{eq:2}
&\tr{\wb_1\xb^{\pm 1}\wb_2 \xb^{\pm 1}\wb_3}=\nonumber\\
\nonumber\\
&-\tr{\wb_1\xb^{\pm 2}\wb_2 \wb_3}-\tr{\wb_1\wb_2 \xb^{\pm 2}\wb_3}+\tr{\wb_1\pol(\xb^{\pm 1}, \wb_2)\wb_3}.
\end{align}
However, by subsequently reducing the words having letters with exponent not $\pm 1$, we may freely eliminate generators of the form $\tr{\wb_1\xb^{\pm 1}\wb_2 \xb^{\pm 1}\wb_3}$. In other words, we may assume that for any word in any generator no letter with the same exponent is ever repeated in the same word.

Letting $\wb_3=\xb$ in Equation \eqref{eq:2} we deduce \begin{align}\label{eq:3} &\tr{\wb_1\xb\wb_2 \xb^2}=\nonumber\\ \nonumber\\&-\tr{\wb_2 \xb\wb_1\xb^2}-\tr{\wb_1\wb_2 \xb^3}+\tr{\wb_1\pol(\xb, \wb_2)\xb}.\end{align}  Then reducing words having letters with exponents not $\pm 1$ we find that we need only one of $\tr{\wb_1\xb\wb_2 \xb^{-1}}$ and $\tr{\wb_2 \xb\wb_1\xb^{-1}}$ to generate $\C[\X_r]$.

Putting these reductions together, we have the following description of the generators of $\C[\X_r]$.

\begin{lem}\label{generatorform}
$\C[\X_r]$ is generated by traces of the form:
\begin{align*}&\tr{\xb_i},\tr{\xb_i^{-1}}, \tr{\xb_i\xb_j}, \tr{\xb_i\xb_j\xb_k}, \tr{\xb_i\xb_j^{-1}},\tr{\xb_i^{-1}\xb_j^{-1}},\\
&\tr{\xb_i\xb_j\xb_k^{-1}},\tr{\xb_i\xb_j\xb_k\xb_l},\tr{\xb_i\xb_j\xb_k\xb_l\xb_m},\tr{\xb_i\xb_j\xb_k\xb_l^{-1}},\\
& \tr{\xb_i\xb_j\xb_k\xb_j^{-1}},\tr{\xb_i\xb_j^{-1}\xb_k^{-1}},
\tr{\xb_i^{-1}\xb_j^{-1}\xb_k^{-1}},\tr{\xb_i\xb_j\xb_k^{-1}\xb_l^{-1}},\\
&\tr{\xb_i\xb_j\xb_k^{-1}\xb_j^{-1}},\tr{\xb_i\xb_j\xb_i^{-1}\xb_j^{-1}}, \tr{\xb_i\xb_j\xb_k\xb_l\xb_m^{-1}},\\ 
&\tr{\xb_i\xb_j\xb_k\xb_l\xb_k^{-1}},\tr{\xb_i\xb_j\xb_k\xb_l\xb_m\xb_n},\end{align*} where
$1\leq i\not=j\not=k\not=l\not=m\not=n \leq r$.
\end{lem}

\begin{proof}
The preceding remarks of this section deserve to be listed (referred to as ``the summary''):

\begin{enumerate}
\item If a word has at least one letter with a negative power, we assume (by a cyclic permutation) one of those letters is the last letter in the word. 
 \item No letter is repeated with the same exponent.
 \item All exponents are $\pm 1$.
 \item All words are of weighted length $6$ or less.
\end{enumerate}

For instance, this tells that the most inverses possible in a word is 3 since the length of a letter with  a negative power in a word is counted twice and the maximal total length is $6$.  

Since exponents are always $\pm 1$, when we say ``a positive exponent'' or ``a negative exponent'' we always mean ``an exponent of $1$'' or ``an exponent of $-1$,'' respectively. 

The description of the possible generator types is organized in tables, and separated by total length and number of negative exponents represented in the word corresponding to the generator.  All letters in a word denote generic unimodular matrices.

We begin with the words of length $4$ or less, since the summary is all that is necessary to describe these generator types without further comment.

\begin{table}[!ht]
\begin{center}
\begin{tabular}{|c|c|c|}
\hline
Length & Negative Exponents & Word Type\\
\hline
$1$ & $0$ & $\tr{\xb}$\\
\hline
$2$ & $0$ & $\tr{\xb\yb}$\\
\hline
$2$ & $1$ & $\tr{\xb^{-1}}$\\
\hline
$3$ & $0$ & $\tr{\xb\yb\zb}$\\
\hline 
$3$ & $1$ & $\tr{\xb\yb^{-1}}$\\
\hline
$4$ & $0$ & $\tr{\wb\xb\yb\zb}$\\
\hline
$4$ & $1$ & $\tr{\xb\yb\zb^{-1}}$\\
\hline
$4$ & $2$ & $\tr{\xb^{-1}\yb^{-1}}$\\
\hline
\end{tabular}
\caption{This table lists the generator types that are in terms of words of length four or less.}\label{4orless}
\end{center}
\end{table}

There is no possibility for the generators listed in Table $\ref{4orless}$ to have any letter coexist in the same word with that letter's inverse (because of cyclic reduction).  With generators of length $5$ and $6$ this becomes both possible and necessary.  

We now address the generators in words of length $5$.

\begin{table}[!ht]
\begin{center}
\begin{tabular}{|c|c|c|}
\hline
Length & Negative Exponents & Word Type\\
\hline
$5$ & $0$ & $\tr{\vb\wb\xb\yb\zb}$\\
\hline
$5$ & $1$ & $\tr{\wb\xb\yb\zb^{-1}}$\\
\hline
$5$ & $1$ & $\tr{\wb\xb\yb\xb^{-1}}$\\
\hline
$5$ & $2$ & $\tr{\xb\yb^{-1}\zb^{-1}}$\\
\hline
\end{tabular}
\caption{This table lists the generator types that are in terms of words of length 5.}\label{length5}
\end{center}
\end{table}

Again, cyclic reduction and the summary together force only the generator types listed in Table $\ref{length5}$ to be sufficient to generate the coordinate ring and have word length $5$.

For instance, if a length $5$ word has two letters with negative exponents it must have exactly one other letter without a negative exponent, since the letters with negative exponents have a weighted length of $2$.  Therefore, the word is in one of the following forms: $\xb^{-1}\yb\zb^{-1}$, $\yb\zb^{-1}\xb^{-1}$, or $\zb^{-1}\xb^{-1}\yb$.  However, all three forms give the same generator type since $$\tr{\xb^{-1}\yb\zb^{-1}}=\tr{\yb\zb^{-1}\xb^{-1}}=\tr{\zb^{-1}\xb^{-1}\yb};$$ exemplifying {\it cyclic equivalence.}

We now address the generators that are in terms of words of length $6$.  This case will require one further reduction formula beyond the summary to complete the table.  We give the table here:

\begin{table}[!ht]
\begin{center}
\begin{tabular}{|c|c|c|}
\hline
Length & Negative Exponents & Word Type\\
\hline
$6$ & $0$ & $\tr{\ub\vb\wb\xb\yb\zb}$\\
\hline
$6$ & $1$ & $\tr{\vb\wb\xb\yb\zb^{-1}}$\\
\hline
$6$ & $1$ & $\tr{\vb\wb\xb\yb\xb^{-1}}$\\
\hline
$6$ & $2$ & $\tr{\wb\xb\yb^{-1}\zb^{-1}}$\\
\hline
$6$ & $2$ & $\tr{\wb\xb\yb^{-1}\xb^{-1}}$\\
\hline
$6$ & $2$ & $\tr{\yb\xb\yb^{-1}\xb^{-1}}$\\
\hline
$6$ & $3$ & $\tr{\xb^{-1}\yb^{-1}\zb^{-1}}$\\
\hline
\end{tabular}
\caption{This table lists the generator types that are in terms of words of length $6$.}\label{length6}
\end{center}
\end{table}

For words having one letter with a negative exponent, obviously there is at most one letter repeated with its inverse.  Cyclic equivalence forces it to be either $\tr{\vb\wb\xb\yb\xb^{-1}}$ or $\tr{\yb\xb\vb\wb\xb^{-1}}$, putting the inverse letter at the end of each respective word.  However, Equation $\eqref{eq:3}$ mandates that we need only one of these.

For words having two letters with negative exponents, once we establish that such letters may be assumed adjacent, the three possibilities listed in Table $\ref{length6}$ are forced to be sufficient to generate the ring after taking into consideration possible cyclic permutations of the listed generator types.

To complete our description of the generators in Table $\ref{length6}$ having two letters with negative exponents, it now remains to consider generators of the form  $\tr{\ub^{-1}\wb\xb^{-1}\zb}$.  Using the algorithm in \cite{L2} (see Appendix A) for reducing traces of words of length 7 or more to those of length 6 or less, one 
computes that $$\tr{\ub\vb\wb\xb\yb\zb}+\tr{\ub\vb\wb\yb\xb\zb}+\tr{\vb\ub\wb\xb\yb\zb}+\tr{\vb\ub\wb\yb\xb\zb}$$ has trace degree $5$; that is, can be expressed as a polynomial in traces of words no longer than 5.  Setting $\ub=\vb$ and $\xb=\yb$ and subsequently interchanging words with squares to those with inverses, we find generators of the form $\tr{\ub^{-1}\wb\xb^{-1}\zb}$ can be freely eliminated; that is, inverses can be assumed to be adjacent.

Clearly, words with three letters all having negative exponents must be in form $\tr{\xb^{-1}\yb^{-1}\zb^{-1}}$.
\end{proof}

\begin{rem}
In \cite{L2} these $19$ generator types are describes plus one additional one that is not necessary.  Namely, Equation \eqref{eq:3} shows us that only one of the generators 
$\tr{\xb_i\xb_j\xb_k\xb_l\xb_k^{-1}}$ and $\tr{\xb_i\xb_k\xb_j\xb_l\xb_k^{-1}}$ is needed.
\end{rem}

\subsection{Second Reduction:  Minimal Generators}\label{secondreduc}

To prove Theorem $\ref{maintheorem}$, it remains to count how many of each type of generator listed Lemma $\ref{generatorform}$ is necessary to generate the ring of invariants $\C[\X_r]$.  We devote the remainder of this section to this goal. 

Using the representation theory of $\mathrm{GL}(r,\C)$, Abeasis and Pittaluga determined in \cite{AP} a method to count the minimal number of generators with respect to word length and with respect to the invariants of \emph{arbitrary} matrices.  Their algorithm is viable for any size matrix, but relies on computer computations that were only implementable for $3\times 3$ matrices.  Additionally, they also derived explicit highest weight vectors which can be used to write down an explicit minimal set of generators for $\glm{3}^{\times r}\aq \G$.  In this case, minimal generators have also been determined by \cite{M,MS1,MS2,MS3} and also by \cite{Lo2,Lo3}.

We now review the method of \cite{AP}.  We then use it to count each type of generator for the unimodular invariants under consideration in this paper.

Recall the notation of Section $\ref{sec2}$:  $$\Y_r=\qt{\glm{3}} \text{ and } \X_r=\qt{\SLm{3}}.$$  As noted earlier $\C[\Y_r]$ is a $\C$-algebra generated by the functions 
$(\xb_1,...,\xb_r)\mapsto \tr{\xb_{i_1}\cdots\xb_{i_k}}$, where $k\leq 6$.  The ring of invariants $\C[\Y_r]=\C[x_{ij}^k]^\G$ is a subring of a connected multigraded ring, and each generator $\tr{\xb_{i_1}\cdots\xb_{i_k}}$ is a homogeneous polynomial of degree $k$ where each monomial has exactly one matrix element from each of the represented generic matrices $\xb_{i_j},$ for $1\leq j\leq k$.  Consequently, $\C[\Y_r]$ is connected and multigraded by the degrees in $\xb_1,...,\xb_r$.   Note that the ring of unimodular invariants $\C[\X_r]$ is not graded.  However, Proposition $\ref{slgl}$ implies that it is filtered since the ideal identifies a homogeneous polynomial of degree $3$ (the determinant) with the degree $0$ polynomial $1$.

The group $\GLm{r}$ acts on $\C[\Y_r]$:  $$g\cdot \tr{\xb_{i_1}\cdots \xb_{i_k}}=\tr{\sum_{i}g_{i_1i}\xb_i\cdots\sum_{i}g_{i_ki}\xb_i},$$ and preserves degree.  Therefore it acts on the positive terms $\C[\Y_r]^+$ and on the vector space $\C[\Y_r]^+/\left(\C[\Y_r]^+\right)^2$.  Consequently, a basis for this vector space pulls back to a set of minimal generators for $\C[\Y_r]$ as a $\C$-algebra.

Abesis and Pittaluga determine the irreducible subspaces of this action on $\C[\Y_r]^+/\left(\C[\Y_r]^+\right)^2$ by highest weight (see \cite{FH} for background in representation theory).  In these terms they construct a minimal basis.  In particular, the set of weights of the irreducible subspaces is $$\{(1),(2),(1^3),(3),(2^2),(2,1^2),(1^5),(3,1^2),(2^2,1),(3^2),(3,1^3)\},$$ where, for instance, when $r=5$ the expression $(1^3)$ denotes the partition $(1,1,1,0,0)$.  Moreover, the sum of the entries in a weight correspond to the degree of the generators, since the weights correspond to Young symmetrizers acting on the generic matrices.  The dimension of these irreducible representations is known classically.  An irreducible $\GLm{r}$ representation having partition $(\lambda_1,...,\lambda_r)$ is $$\prod_{1\leq i<j\leq r}\frac{\lambda_i-\lambda_j+j-i}{j-i}.$$  Note that with respect to the language used in \cite{FH} these are the \emph{conjugate} partitions of the ones used in \cite{AP}.  Naturally, if $r$ is less than the length of the partition then the basis is empty.

We will determine explicit reductions for each generator type and use the dimension of the irreducible representations as a means to count when we have enough relations.

Proposition $\ref{numbgen}$ and the work of \cite{AP} together will allow us to determine the number of minimal generators in $\C[\X_r]$.  However, we wish to know more.  We will further write down explicit (at times algorithmic) reductions to take the set of generators of a certain type to a minimal sufficient subset of generators of that type.

We use the same notation, $\N{x,y}$, to denote the number of generators with respect to the free group of rank $r$ of word length $x$ in $y$ letters in $\C[\X_r]$ (as opposed to $\C[\Y_r]$).  Recall that any letter with an exponent of $-1$ has length $2$.

Again to organize the information, as it is copious, we use tables.  We proceed with the easiest and most immediate cases.  When $r=1$, $\SL\aq\SL\approx \C^2$ is affine and thus
$\C[\SL\aq \SL]\approx \C[\tr{\xb},\tr{\xb^{-1}}].$  Consequently, Proposition $\ref{count}$ gives the following table:

\begin{table}[!ht]
\begin{center}
\begin{tabular}{|c|c|}
\hline
Minimal Number & Generators \\
\hline
$\N{1,1}=r$ & $\tr{\xb_1},..,\tr{\xb_r}$ \\
\hline
$\N{2,1}=r$ & $\tr{\xb_1^{-1}},...,\tr{\xb^{-1}_r}$\\
\hline
\end{tabular}
\end{center}
\caption{This table lists the minimal generators in words with only one letter.}\label{min1}
\end{table}

There are no further reductions necessary, as all the listed generators are necessary to generate the invariant ring.  Note that this corroborates with the computation of the dimensions of the irreducible representations.  For in this case the only weight vectors are $(1),(2),$ and $(3)$.  The dimension in each case is trivially $1$ and the basis is $\{\tr{\xb},\tr{\xb^2},\tr{\xb^3}\}$, which projects to $\{\tr{\xb},\tr{\xb^{-1}}\}$.

For $r=2$ we refer to \cite{L2}.  Then $\C[\SL^{\times 2}\aq \SL]$ is minimally generated by:  \begin{align*}
&\{\tr{\xb_1}, \tr{\xb_1^{-1}}, \tr{\xb_2},\tr{\xb_2^{-1}}, \tr{\xb_1\xb_2},\tr{\xb_1\xb_2^{-1}},\\ &\tr{\xb_2\xb_1^{-1}},\tr{\xb_1^{-1}\xb_2^{-1}}, 
\tr{\xb_1\xb_2\xb_1^{-1}\xb_2^{-1}}\}.\end{align*}
The corresponding minimal generators for $\C[\glm{3}^{\times 2}\aq\SL]$ are:
\begin{align*}
&\{\tr{\xb_1}, \tr{\xb_1^{2}}, \tr{\xb_1^3}, \tr{\xb_2},\tr{\xb_2^{2}}, \tr{\xb^3}, \tr{\xb_1\xb_2},\\ 
&\tr{\xb_1\xb_2^{2}}, 
\tr{\xb_2\xb_1^{2}},\tr{\xb_1^{2}\xb_2^{2}}, \tr{\xb_1\xb_2\xb_1^{2}\xb_2^{2}}\}.
\end{align*}
The weight vectors for the $r=2$ case are $(1,0),(2,0),(3,0),(2,2)$ and $(3,3)$.  Their dimensions are respectively computed to be $2,3,4,1,1$ which add to $11$; the number of minimal generators for this case.  From the $r=1$ case we can account for the first $6$ of these, and the rest are essential to $r=2$.  Moreover, from Proposition $\ref{numbgen}$ we would expect $11-2=9$ generators for $\C[\X_2]$; the number listed. 

We then have the following table:

\begin{table}[!ht]
\begin{center}
\begin{tabular}{|c|c|}
\hline
Minimal Number & Generators\\
\hline
$\N{2,2}=\binom{r}{2}$ & $\{\tr{\xb_i\xb_j}\}$, where $1\leq i<j \leq r$\\
\hline
$\N{3,2}=2\binom{r}{2}$ & $\{\tr{\xb_i\xb_j^{-1}}\}$, where $i\not=j$ and $1\leq i,j \leq r$\\
\hline
$\N{4,2}=\binom{r}{2}$ & $\{\tr{\xb_i^{-1}\xb_j^{-1}}\}$, $1\leq i<j \leq r$\\
\hline
$\N{6,2}=\binom{r}{2}$ & $\{\tr{\xb_i\xb_j\xb_i^{-1}\xb_j^{-1}}\}$, $1\leq i<j \leq r$\\
\hline
\end{tabular}
\end{center}
\caption{This table lists the minimal generators in words with only two letters.}\label{min2}
\end{table}

At this point we make a couple observations.  Rows $1$, $2$, and $3$ above are determined by cyclic reduction alone; all generators of that type are included up to cyclic equivalence.  We are making a choice however, we are choosing for instance $\tr{\xb_1\xb_2}$ over $\tr{\xb_2\xb_1}$; but they are identical.  There is also a choice involved in row $4$, but this one is not trivial.  In \cite{L2} the following relation in $\C[\X_2]$ is derived:

\begin{align}
&\tr{\xb_2\xb_1\xb_2^{-1}\xb_1^{-1}}=\nonumber \\
&\nonumber\\
& -\tr{\xb_1\xb_2\xb_1^{-1}\xb_2^{-1}}+\tr{\xb_1}\tr{\xb_1^{-1}}\tr{\xb_2}\tr{\xb_2^{-1}}
+\tr{\xb_1}\tr{\xb_1^{-1}}\nonumber\\
&+\tr{\xb_2}\tr{\xb_2^{-1}}+\tr{\xb_1\xb_2}\tr{\xb_1^{-1}\xb_2^{-1}}
+\tr{\xb_1\xb_2^{-1}}\tr{\xb_1^{-1}\xb_2}\nonumber\\
&-\tr{\xb_1^{-1}}\tr{\xb_2}\tr{\xb_1\xb_2^{-1}}
-\tr{\xb_1}\tr{\xb_2^{-1}}\tr{\xb_1^{-1}\xb_2}\nonumber\\
&-\tr{\xb_1}\tr{\xb_2}\tr{\xb_1^{-1}\xb_2^{-1}}
-\tr{\xb_1\xb_2}\tr{\xb_1^{-1}}\tr{\xb_2^{-1}}-3.\label{rank2sum}
\end{align}

We note that Equation $\eqref{rank2sum}$ corresponds to a relation between $\tr{\xb\yb\xb^2\yb^2}$ and $\tr{\yb\xb\yb^2\xb^2}$ in $\C[\Y_2]$.

Up to cyclic equivalence there are only three words of length $6$ in two letters, the two just mentioned and $\tr{\xb^3\yb^3}$; the latter most being reducible since it has letters with exponents greater than $2$.  The maximum exponent allowed on any letter in any word even in $\C[\Y_r]$ is $2$, except for $\tr{\xb^3}$ itself.

Hence using Equation $\eqref{rank2sum}$ we may choose the generators in row $4$ of Table $\ref{min2}$.  Therefore, we have accounted for minimality with explicit reductions in the cases $r=1,2$.

We have to consider all cases up to $r=6$ since the maximum word length necessary is $6$, and Proposition $\ref{count}$ shows that both the generator types and the combinatorics are determined in general by the cases $r=1,...,6$ alone.

The highest weight vectors for $r=3$ are:  $\Lambda_3=$ $$\{(1,0,0),(2,0,0),(3,0,0),(2,2,0),(3,3,0),(1,1,1),(2,1,1),(3,1,1),(2,2,1)\}.$$  
For ease of reading, we will often abbreviate ``irreducible representation'' with ``irrep'' from this point on.

Letting $\Y_3=\glm{3}^{\times 3}\aq \G$ and $V_3{(\lambda)}$ be the $\GLm{3}$ irrep of $$V_3=\C[\Y_3]^+/\left(\C[\Y_3]^+\right)^2$$ associated to $\lambda\in \Lambda_3$, we may write $$V_3=\bigoplus_{\lambda\in \Lambda_3}V_3{(\lambda)}.$$

The dimension of each irrep listed in $\Lambda_3$ is, respectively:  $3,6,10,6,10,1,3,6,3$; adding to $48$.  From the previous cases we can account for $21$ generators that form the bases of these irreps plus the $3$ removed generators of the form $\tr{\xb^3}$.  We have exactly $48-21-3=24$ left to find.

The following table categorizes the number and type of these remaining generators (we suppress the indexing for clarity):

\begin{table}[!ht]
\begin{center}
\begin{tabular}{|c|c|}
\hline
Minimal Number & Generators\\
\hline
$\N{3,3}=2\binom{r}{3}$ & $\tr{\xb\yb\zb},\tr{\yb\xb\zb}$\\
\hline
$\N{4,3}=6\binom{r}{3}$ & $\tr{\xb\yb\zb^{-1}}, \tr{\yb\xb\zb^{-1}}$\\
\hline
$\N{5,3}=9\binom{r}{3}$ & $3\binom{r}{3}$: $\tr{\xb\yb\zb\yb^{-1}}$, $6\binom{r}{3}$: $\tr{\xb\yb^{-1}\zb^{-1}}$\\
\hline
$\N{6,3}=7\binom{r}{3}$ & $6\binom{r}{3}$: $\tr{\xb\yb\zb^{-1}\yb^{-1}}$, $\binom{r}{3}$: $\tr{\xb^{-1}\yb^{-1}\zb^{-1}}$\\
\hline
\end{tabular}
\end{center}
\caption{This table lists the minimal generators in words with only three letters.}\label{min3}
\end{table}

The first row of Table $\ref{min3}$ concerns words of length $3$ in $3$ letters.  The words of length $3$ that are accounted for by the cases $r=1$ and $r=2$ number $9$, meaning there are exactly $2$ remaining since there is a total of $\dim V_3(3,0,0)+\dim V_3(1,1,1)=10+1=11$.  However, up to cyclic equivalence there are only $2$ possible words of length $3$ in three letters; both are listed in row $1$.

Row $2$ of Table $\ref{min3}$ concerns words of length $4$ in $3$ letters.  There are $9$ generators corresponding to such words.  The words of length $4$ that are accounted for number $3$, so there are $6$ remaining.  Again, up to cyclic equivalence there are exactly $6$ possibilities; justifying row $2$.

Generators in terms of words of length $5$ first show up in this case so the $9$ claimed to be required need to be described.  Cyclic equivalence gives us exactly $6$ of the generators of type $\tr{\xb\yb^{-1}\zb^{-1}}$.  However, cyclic equivalence alone gives us $6$ more of type $\tr{\xb\yb\zb\yb^{-1}}$. Equation $\eqref{eq:3}$ mandates that we only need half of these.  In other words, for every choice of $\yb$ we can choose either $\tr{\xb\yb\zb\yb^{-1}}$ or $\tr{\zb\yb\xb\yb^{-1}}$.  This gives us the required $9$ and minimality promises there are no further reductions (as it goes with every case).

We come to the fourth row.  There are $3$ generators accounted for in terms of length $6$ words from the two cases when $r$ is in $\{1,2\}$; namely, $$\{\tr{\xb_1\xb_2\xb_1^{-1}\xb_2^{-1}}, \tr{\xb_1\xb_3\xb_1^{-1}\xb_3^{-1}}, \tr{\xb_2\xb_3\xb_2^{-1}\xb_3^{-1}} \}.$$  With a total of $10$ required by the dimension of the irrep with weight $(3,3,0)$; there are $7$ left to find.  Equation $\eqref{eq:3}$ tells us that we can assume that inverses are together, and with that noted, cyclic equivalence alone gives us $6$ generators of type $\tr{\xb\yb\zb^{-1}\yb^{-1}}$.  The last possible generator type is $\tr{\xb^{-1}\yb^{-1}\zb^{-1}}$.  Cyclically, there are $2$ of these.

Using the algorithm in Appendix A for the expression $$\tr{\ub\vb\wb\xb\yb\zb}+\tr{\ub\vb\wb\yb\xb\zb}+\tr{\vb\ub\wb\xb\yb\zb}+\tr{\vb\ub\wb\yb\xb\zb}$$ and identifying $\ub=\vb$, 
$\wb=\xb$ and $\yb=\zb$ we come to an expression of trace degree five:  $2\tr{\ub^2\xb^2\yb^2}+2\tr{\ub^2(\xb\yb)^2}$.  Iteratively replacing squares with inverses reduces this to an expression having 
the same trace degree:  $$2\tr{\ub^{-1}\xb^{-1}\yb^{-1}}+2\tr{\ub^{-1}(\xb\yb)^{-1}}=2\tr{\ub^{-1}\xb^{-1}\yb^{-1}}+2\tr{\ub^{-1}\yb^{-1}\xb^{-1}},$$ which shows that we can choose the order of any 
three such letters and so finishes the proof that there are $\binom{r}{3}$ generators of the form $\tr{\xb^{-1}\yb^{-1}\zb^{-1}}$.

We now move on the generators in terms of four letters.  The table is as follows:

\begin{table}[!ht]
\begin{center}
\begin{tabular}{|c|c|}
\hline
Minimal Number & Generators\\
\hline
$\N{4,4}=5\binom{r}{4}$ & $\tr{\wb\xb\yb\zb},\tr{\wb\xb\zb\yb},\tr{\wb\yb\xb\zb},$\\
& $\tr{\wb\yb\zb\xb},\tr{\wb\zb\xb\yb}$\\
\hline
$\N{5,4}=20\binom{r}{4}$ & $\tr{\wb\xb\yb\zb^{-1}}$\\
\hline
$\N{6,4}=26\binom{r}{4}$ & $18\binom{r}{4}$: $\tr{\wb\xb\yb^{-1}\zb^{-1}}$, $8\binom{r}{4}$: $\tr{\wb\xb\yb\zb\yb^{-1}}$\\
\hline
\end{tabular}
\end{center}
\caption{This table lists the minimal generators in words with only four letters.}\label{min4}
\end{table}

The irreducible spaces consisting of degree four generators have partitions: $(2,2,0,0)$ and $(2,1,1,0)$; having total dimension $35$.  From $r=2$ we have $\binom{4}{2}=6$ generators in the form $\tr{\xb^{-1}\yb^{-1}}$; $r=3$ gives an additional $6\binom{4}{3}=24$ of type $\tr{\xb^{-1}\yb\zb}$.  This leaves exactly $5$ generators to find.  There is only one further generator type of length four in four letters:  $\tr{\wb\xb\yb\zb}$.  Cyclically there are $6$ of these.

We want to explicitly construct the minimal generating set, so we want the relation.  In \cite{L2} the following relation is derived:
\begin{align*}
&\xb\zb\yb+\zb\xb\yb+\yb\xb\zb+\yb\zb\xb+\xb\yb\zb+\zb\yb\xb=\\ \\
&\pol(\xb+\zb,\yb)-\pol(\xb,\yb)-\pol(\zb,\yb);\end{align*}  it is the full polarization of the Cayley-Hamilton equation.
Thus,
\begin{align}\label{cyclicsum}
&\tr{\wb\xb\zb\yb}+\tr{\wb\zb\xb\yb}+\tr{\wb\yb\xb\zb}+\nonumber\\
&\tr{\wb\yb\zb\xb}+\tr{\wb\xb\yb\zb}+\tr{\wb\zb\yb\xb}=\nonumber\\ \nonumber\\
&\tr{\wb\left(\pol(\xb+\zb,\yb)-\pol(\xb,\yb)-\pol(\zb,\yb) \right)}, 
\end{align}
which allows us to eliminate exactly one of the six generators on the left-hand-side of Equation $\eqref{cyclicsum}$.  This validates the content of row 1 of Table $\ref{min4}$.

We note for later use that Equation $\ref{cyclicsum}$, with respect to general substitutions of words for $\wb,\xb,\yb,$ and $\zb$, has trace degree (recall the trace degree is the largest word length in a trace expression) at most the length of $\wb$ plus one less than the sum of the lengths of $\xb,\yb,\zb$.  This fact is apparent by inspection of Equation $\eqref{eq:1}$, and otherwise expressed $$\left|\tr{\wb\sum\xb\yb\zb}\right|\leq |\wb|+\left(|\xb|+|\yb|+|\zb|-1  \right).$$

With respect to row $2$ in Table $\ref{min4}$, there are $56$ generators of length $5$ coming from the weights $(2,2,1)$ and $(3,1,1)$.  From the cases $r=1,2,3$ we account for $9\binom{4}{3}=36$ of these, leaving $20$ to find.
The only generator type in four letters of length $5$ is $\tr{\wb\xb\yb\zb^{-1}}$.  Cyclically there are $24$ possibilities.  Replacing $\wb$ with $\wb^{-1}$ in Equation $\eqref{cyclicsum}$ gives row $2$, since this provides four reductions (one for each choice of the last letter to have the negative exponent).

Row $3$ in Table $\ref{min4}$ describes the generators of length $6$ in four letters.  The corresponding weights are $(3,3)$ and $(3,1,1,1)$ giving a total of $60$ generators.  The previous cases account for $34$ of these, leaving $26$ to derive.  Lemma \ref{generatorform} shows that the only generator types of length $6$ in four letters are $\tr{\wb\xb\yb^{-1}\zb^{-1}}$ and $\tr{\wb\xb\yb\zb\yb^{-1}}$.

We first address the generators in the form $\tr{\wb\xb\yb\zb\yb^{-1}}$.  The relations between generators of this type and the generators that have a form given by permutating the letters of $\tr{\wb\xb\yb\zb\yb^{-1}}$ are given in the proof and preceding remarks of Lemma \ref{generatorform}.  We count the possible number of generators of this type in this form only.  In other words, we do not again explain why we do not need $\tr{\zb\yb\wb\xb\yb^{-1}}$.

There are four choices for the letter $\yb$ and there are three remaining for the letter $\zb$.  We first show that we can choose an order for $\wb\xb$ and demonstrate the relation between differing orders.  Then for each choice of $\yb$ we show there is a sum relation among the remaining three choices for $\zb$ which gives four further reductions.

Using the algorithm in Appendix A for the expression $$\tr{\ub\vb\wb\xb\yb\zb}+\tr{\ub\vb\wb\yb\xb\zb}+\tr{\vb\ub\wb\xb\yb\zb}+\tr{\vb\ub\wb\yb\xb\zb}$$ and setting $\xb=\wb=\zb$ gives the expression $$\tr{\ub\vb\xb^2\yb\xb}+\tr{\ub\vb\xb\yb\xb^2}+\tr{\vb\ub\xb^2\yb\xb}+\tr{\vb\ub\xb\yb\xb^2},$$ having trace degree $5$.  Exchanging squares for inverses and cyclically permuting the letters gives that $$\tr{\yb\xb\ub\vb\xb^{-1}}+\tr{\ub\vb\xb\yb\xb^{-1}}+\tr{\yb\xb\vb\ub\xb^{-1}}+\tr{\vb\ub\xb\yb\xb^{-1}}$$ has trace degree $5$ as well.  Then using Equation $\eqref{eq:3}$ we 
reduce this expression to the trace degree $5$ expression 
\begin{align}
&2\tr{\ub\vb\xb\yb\xb^{-1}}+2\tr{\vb\ub\xb\yb\xb^{-1}},\label{6in4orderrel}
\end{align}
which is what we wished to derive.  In other words we can always choose an order for the first two letters of such a generator and thus there are no more than $12$ such generators necessary.

Now, we show there is a sum relation giving four further reductions, as required to show there is no more than $8$ generators of this type necessary.

Setting $\wb=\xb^{-1}\ub$ in Equation $\eqref{cyclicsum}$ and cyclically permuting letters gives a relation of trace degree $5$:
\begin{align*}
&\tr{\ub\xb\zb\yb\xb^{-1}}+\tr{\ub\zb\xb\yb\xb^{-1}}+\tr{\ub\yb\xb\zb\xb^{-1}}\\
&+\tr{\ub\yb\zb}+\tr{\ub\xb\yb\zb\xb^{-1}}+\tr{\ub\zb\yb}.
\end{align*}
By subtracting generators with words of length three and using Equation $\eqref{eq:3}$ we come to a relation of like trace degree:  \begin{equation}\label{eq:step}
\tr{\zb\yb\xb\ub\xb^{-1}}+\tr{\ub\zb\xb\yb\xb^{-1}}+\tr{\ub\yb\xb\zb\xb^{-1}}+\tr{\yb\zb\xb\ub\xb^{-1}}.
\end{equation}
However, we just showed (see Expression $\eqref{6in4orderrel}$) that expressions of the form $$2\tr{\ub\vb\xb\yb\xb^{-1}}+2\tr{\vb\ub\xb\yb\xb^{-1}},$$ when both terms are taken together, can be freely eliminated.  Hence, the outer two summands in Expression $\ref{eq:step}$ can be eliminated since they are equal to an expression in terms of generators of smaller word length.  We come to the conclusion that $$\tr{\ub\zb\xb\yb\xb^{-1}}+\tr{\ub\yb\xb\zb\xb^{-1}}$$ has trace degree $5$ as well. Therefore, for every group of three generators $$\{\tr{\ub\zb\xb\yb\xb^{-1}},\tr{\zb\yb\xb\ub\xb^{-1}},\tr{\ub\yb\xb\zb\xb^{-1}} \}$$ there is one reduction, and there are four such groups among the twelve remaining 
generators of this type.  Consequently, we have established four further reductions and thus there are no more than $8$ generators of this type necessary.

The other possible generator form not addressed yet from row $3$ in Table $\ref{min4}$
is $\tr{\wb\xb\yb^{-1}\zb^{-1}}$.  We must show there are no more than $18$ necessary.  Then the $8+18$ generators so derived will give the minimal number of $26$, which will imply that no further generators can be eliminated from the sets of either type.

There are exactly $24$ possible generators in the form $\tr{\wb\xb\yb^{-1}\zb^{-1}}$.  We must find $6$ reductions.

Again using the algorithm in Appendix A for the expression $$\tr{\ub\vb\wb\xb\yb\zb}+\tr{\ub\vb\wb\yb\xb\zb}+\tr{\vb\ub\wb\xb\yb\zb}+\tr{\vb\ub\wb\yb\xb\zb}$$ but this time setting $\xb=\yb$ 
and $\zb=\wb$ gives that the expression $$2\tr{\ub\vb\wb\xb^2\wb}+2\tr{\vb\ub\wb\xb^2\wb}$$ is equal to a polynomial in generators of word length at most $5$.  Using Equation $\eqref{eq:2}$ reduces this to $$2\tr{\ub\vb\wb^2\xb^2}+2\tr{\vb\ub\wb^2\xb^2}+2\tr{\ub\vb\xb^2\wb^2}+2\tr{\vb\ub\xb^2\wb^2};$$  trading squares for inverses provides $$2\tr{\ub\vb\wb^{-1}\xb^{-1}}+2\tr{\vb\ub\wb^{-1}\xb^{-1}}+2\tr{\ub\vb\xb^{-1}\wb^{-1}}+2\tr{\vb\ub\xb^{-1}\wb^{-1}}$$ also has trace degree $5$.  Hence among every 
group of four generators $$\{\tr{\ub\vb\wb^{-1}\xb^{-1}}, \tr{\vb\ub\wb^{-1}\xb^{-1}}, \tr{\ub\vb\xb^{-1}\wb^{-1}}, \tr{\vb\ub\xb^{-1}\wb^{-1}}\}$$ we obtain exactly one relation allowing for six reductions since for any four letters there are $6$ such collections (count multidegree).  And so we have the required $24-6=18$ generators of this type.

Next, we address generators in five letters.  From Lemma \ref{generatorform} the only generator types in five letters that have not been accounted for by the cases $r=1,2,3,4$ are $\tr{\ub\vb\wb\xb\yb}$ of length five and $\tr{\vb\wb\xb\yb\zb^{-1}}$ of length six.  The table follows:

\begin{table}[!ht]
\begin{center}
\begin{tabular}{|c|c|}
\hline
Minimal Number & Generators\\
\hline
$\N{5,5}=12\binom{r}{5}$ & $\tr{\ub\vb\wb\xb\yb}$\\
\hline
$\N{6,5}=35\binom{r}{5}$ & $\tr{\vb\wb\xb\yb\zb^{-1}}$\\
\hline
\end{tabular}
\end{center}
\caption{This table lists the minimal generators in words with only five letters.}\label{min5}
\end{table}

We first address row $1$ of Table $\ref{min5}$.  The corresponding weights for the irreps corresponding to length $5$ generators are $(1^5), (3,1^2), (2^2,1)$.  The total dimension of these three irreps is $202$.  There are $9\binom{5}{3}+20\binom{5}{4}=190$ such generators accounted for by our previous work (the cases $r=1,2,3,4$).  This leaves $12$ generators of type $\tr{\ub\vb\wb\xb\yb}$ out of $24$ possible after considering cyclic permutations.  Thus we need $12$ reductions.

We now introduce some new notation to express longer formulas with less symbols.  Let $$\T(1,2,3)=\tr{\xb_1\xb_2\xb_3},\ \T(-1,3)=\tr{\xb_1^{-1}\xb_3},\ \text{etc.}$$  Also, let $\sum \xb\yb\zb$ denote $\xb\zb\yb+\zb\xb\yb+\yb\xb\zb+\yb\zb\xb+\xb\yb\zb+\zb\yb\xb$.

Cyclically we can assume that all generators of type $\tr{\ub\vb\wb\xb\yb}$ in the letters $\{\xb_1,...,\xb_5\}$ are in the form $\T(1,i,j,k,l)$; that is, we assume that $\xb_1$ is the first letter of the word.  This choice of cyclically permuting the letters determines $24$ generators of this type out of the total possible $5!=120$.

We begin addressing the contents of Table $\ref{min5}$ by deriving some useful formulas.  Using the fundamental relation from Appendix A for
\begin{align*}
&3\tr{\xb_1\xb_2\xb_3\xb_4\xb_5\xb_6}+3\tr{\xb_1\xb_2\xb_3\xb_5\xb_4\xb_6}\\
&+3\tr{\xb_2\xb_1\xb_3\xb_4\xb_5\xb_6}+3\tr{\xb_2\xb_1\xb_3\xb_5\xb_4\xb_6}
\end{align*}
and setting $\xb_6=\id$ we derive (using {\it Mathematica}):
\begin{align*}
&3 \T(1, 2, 3, 4, 5) + 3 \T(1, 2, 3, 5, 4) + 3 \T(2, 1, 3, 4, 5) + 3 \T(2, 1, 3, 5, 4)=\\ \\
&3 \T(1) \T(2) \T(5) \T(3, 4) - 3 \T(5) \T(1, 2) \T(3, 4) - 3 \T(2) \T(1, 5) \T(3, 4)\\
&- 3 \T(1) \T(2, 5) \T(3, 4) - 3 \T(2) \T(5) \T(1, 3, 4) + 3 \T(2, 5) \T(1, 3, 4) \\
&- 3 \T(1) \T(5) \T(2, 3, 4)+ 3 \T(1, 5) \T(2, 3, 4) - 4 \T(1, 2) \T(3, 4, 5) \\
&+ 3 \T(3, 4) \T(5, 1, 2) + 3 \T(3, 4) \T(5, 2, 1) - 3 \T(1) \T(2) \T(5, 3, 4) \\
&+ 7 \T(1, 2) \T(5, 3, 4)+ 3 \T(5) \T(1, 2, 3, 4) - 2 \T(2) \T(1, 3, 4, 5) \\
&+ 3 \T(2) \T(1, 5, 3, 4) + 3 \T(5) \T(2, 1, 3, 4) - 2 \T(1) \T(2, 3, 4, 5) \\
&+ 3 \T(1) \T(2, 5, 3, 4) + 5 \T(2) \T(5, 1, 3, 4) + 5 \T(1) \T(5, 2, 3, 4) \\
&- 3 \T(1, 2, 3, 4, 5) + 3 \T(1, 2, 3, 5, 4) - 3 \T(1, 2, 5, 3, 4)- 3 \T(1, 5, 2, 3, 4)\\
& - 3 \T(2, 1, 3, 4, 5) + 3 \T(2, 1, 3, 5, 4)- 3 \T(2, 1, 5, 3, 4)- 3 \T(2, 5, 1, 3, 4)\\
& + 3 \T(5, 1, 2, 3, 4) + 3 \T(5, 2, 1, 3, 4).
\end{align*}

Then bringing all length $5$ generators of this relation to the left side yields:

\begin{align}\label{fundfive}
&3 \T(1, 2, 3, 4, 5) + 3 \T(1, 2, 5, 3, 4) + 3 \T(1, 5, 2, 3, 4) \nonumber\\ &+ 3 \T(1, 3, 4, 5, 2) + 3 \T(1, 5, 3, 4, 2) +3 \T(1, 3, 4, 2, 5) =\nonumber\\ \nonumber\\
&3 \T(1) \T(2) \T(5) \T(3, 4) - 3 \T(5) \T(1, 2) \T(3, 4) - 3 \T(2) \T(1, 5) \T(3, 4)\nonumber\\
& - 3 \T(1) \T(2, 5) \T(3, 4) - 3 \T(2) \T(5) \T(1, 3, 4) + 3 \T(2, 5) \T(1, 3, 4) \nonumber\\
&- 3 \T(1) \T(5) \T(2, 3, 4) + 3 \T(1, 5) \T(2, 3, 4)- 4 \T(1, 2) \T(3, 4, 5) \nonumber\\
&+ 3 \T(3, 4) \T(5, 1, 2) + 3 \T(3, 4) \T(5, 2, 1) - 3 \T(1) \T(2) \T(5, 3, 4) \nonumber\\
&+ 7 \T(1, 2) \T(5, 3, 4)+ 3 \T(5) \T(1, 2, 3, 4) - 2 \T(2) \T(1, 3, 4, 5) \nonumber\\
&+ 3 \T(2) \T(1, 5, 3, 4) + 3 \T(5) \T(2, 1, 3, 4) - 2 \T(1) \T(2, 3, 4, 5) \nonumber \\
&+ 3 \T(1) \T(2, 5, 3, 4) + 5 \T(2) \T(5, 1, 3, 4) + 5 \T(1) \T(5, 2, 3, 4).
\end{align}
Applying the following permutations to Equation $\eqref{fundfive}$ we derive $8$ further relations:  $\{(1),(34),(23),(24),(45),(234),(243),(345)\}$.  We will see that they are independent of each other.  For now, we need four more relations in these generators.

Using Equation $\eqref{cyclicsum}$ and letting $\wb=\ub\vb$ we have a relation for the expression $\tr{\ub\vb\sum \xb\yb\zb}=$ 
\begin{align*}
&\tr{\ub\vb\xb\zb\yb}+\tr{\ub\vb\zb\xb\yb}+\tr{\ub\vb\yb\xb\zb}+\\
&\tr{\ub\vb\yb\zb\xb}+\tr{\ub\vb\xb\yb\zb}+\tr{\ub\vb\zb\yb\xb}.
\end{align*}

There are four relations of the form $\tr{\xb_1\xb_j\sum \xb_k\xb_l\xb_m}$ where $$2\leq 
k,l,m\leq 5, \text{ and } k\not= l\not= m\not= j,$$ corresponding to the four choices for $j=2,3,4,5$.

These additional four relations provide a total, with the other eight from Equation $\eqref{fundfive}$, of twelve relations in $24$ variables.

Since the relations are linear in these variables on the left-hand-side and in all cases the right-hand-side is of strictly smaller trace degree, we form a $12 \times 24$ matrix where all entries are either $1$ or $0$ (we divide Equation $\eqref{fundfive}$ by $3$).  The rank of this matrix is computed in {\it Mathematica} to be $12$.  Thus we have established that these $12$ relations are independent and we can make exactly $12$ reductions leaving us with the required $12$ generators of this type.  The complement of any set of twelve pivot columns of this matrix gives a minimal set of these generators.

The $12\times 24$ matrix is as follows:

{\tiny
$$
\left(
\begin{array}{llllllllllllllllllllllll}
 1 & 0 & 1 & 0 & 0 & 0 & 1 & 0 & 1 & 0 & 0 & 0 & 1 & 0 & 1 & 0 & 0 & 0 & 0 & 0 & 0 & 0 & 0 & 0 \\
 0 & 0 & 0 & 0 & 1 & 1 & 0 & 1 & 1 & 0 & 0 & 0 & 0 & 1 & 1 & 0 & 0 & 0 & 0 & 0 & 0 & 0 & 0 & 0 \\
 0 & 1 & 1 & 0 & 0 & 0 & 0 & 0 & 0 & 0 & 1 & 1 & 1 & 0 & 0 & 1 & 0 & 0 & 0 & 0 & 0 & 0 & 0 & 0 \\
 1 & 0 & 0 & 1 & 0 & 0 & 1 & 0 & 0 & 1 & 0 & 0 & 0 & 0 & 0 & 0 & 1 & 1 & 0 & 0 & 0 & 0 & 0 & 0 \\
 0 & 1 & 0 & 0 & 1 & 0 & 0 & 1 & 0 & 0 & 1 & 0 & 0 & 0 & 0 & 0 & 0 & 0 & 1 & 0 & 1 & 0 & 0 & 0 \\
 0 & 0 & 1 & 1 & 0 & 0 & 1 & 0 & 0 & 0 & 1 & 0 & 0 & 0 & 0 & 0 & 0 & 0 & 0 & 1 & 1 & 0 & 0 & 0 \\
 1 & 0 & 0 & 0 & 1 & 0 & 0 & 0 & 1 & 1 & 0 & 0 & 0 & 0 & 0 & 0 & 0 & 0 & 1 & 0 & 0 & 1 & 0 & 0 \\
 0 & 0 & 0 & 1 & 1 & 0 & 0 & 0 & 0 & 0 & 0 & 0 & 0 & 1 & 0 & 0 & 0 & 1 & 0 & 0 & 1 & 1 & 0 & 0 \\
 1 & 1 & 0 & 0 & 0 & 0 & 0 & 0 & 0 & 0 & 0 & 0 & 1 & 0 & 0 & 0 & 1 & 0 & 1 & 0 & 0 & 0 & 1 & 0 \\
 0 & 0 & 1 & 0 & 0 & 1 & 0 & 0 & 0 & 0 & 0 & 0 & 0 & 0 & 1 & 1 & 0 & 0 & 0 & 1 & 0 & 0 & 0 & 1 \\
 0 & 0 & 0 & 0 & 0 & 0 & 0 & 0 & 0 & 1 & 0 & 1 & 0 & 0 & 0 & 1 & 0 & 1 & 0 & 0 & 0 & 1 & 0 & 1 \\
 0 & 1 & 0 & 0 & 0 & 1 & 0 & 1 & 0 & 0 & 0 & 1 & 0 & 0 & 0 & 0 & 0 & 0 & 0 & 0 & 0 & 0 & 1 & 1
\end{array}
\right)$$}

One such solution set, from the complement of the pivot set, is:
\begin{center}
$\{\T(1, 3, 2, 5, 4), \T(1, 3, 5, 4, 2), \T(1, 4, 2, 5, 3), \T(1, 4, 3, 2, 5),$ \\
$  \T(1, 4, 3, 5, 2), \T(1, 4, 5, 2, 3), \T(1, 4, 5, 3, 2), \T(1, 5, 2, 4, 3),$ \\
$  \T(1, 5, 3, 2, 4), \T(1, 5, 3, 4, 2), \T(1, 5, 4, 2, 3), \T(1, 5, 4, 3, 2)\}.$
\end{center}

We now turn to row $2$ in Table $\ref{min5}$; listing the minimal generators in $5$ letters of length $6$.  The only generators of length $6$ in five letters are the generators of type $\tr{\vb\wb\xb\yb\zb^{-1}}$.  The weight vectors for the irreps corresponding to the length six generators are $\{(3,1,1,1),(3,3)\}$.  The sum of the dimesions of these two irreps totals $245$.  The cases $r=1,2,3,4$ account for $26\binom{5}{4} + 7\binom{5}{3} + \binom{5}{2}=210$ of these generators, leaving exactly $35$ generators of this type to describe.

There are five choices for the letter with the negative exponent in the generator $\tr{\vb\wb\xb\yb\zb^{-1}}$, and we can assume it is always the last letter of the word.  For every such choice there are $24$ choices for the other four letters.  It suffices to show there are $17$ independent relations for every one of the five choices for $\zb^{-1}$.  Then we will be left with $(24-17=7)\times 5=35$ generators, as required.

First, we provide $12$ easy reductions; that is, expressions whose trace degree is less than the largest word length in the expression.  Such an expression, when uniformly in terms of generators of a fixed type, permits one to eliminate one of the generators in the expression.  Indeed, letting $\xb_1=\xb_2$ in the fundamental relation (see Appendix A) gives a reduction formula (recall that this means it is equal to an expression having lesser trace degree) for the expression 
\begin{align*}
&\tr{\xb_1^2\xb_3\xb_4\xb_5\xb_6}+\tr{\xb_1^2\xb_3\xb_5\xb_4\xb_6}\\
&+\tr{\xb_1^2\xb_3\xb_4\xb_5\xb_6}+\tr{\xb_1^2\xb_3\xb_5\xb_4\xb_6}.
\end{align*}
Switching squares for inverses, cyclically permuting and re-indexing gives a further reduction for \begin{equation}
 2\tr{\xb_1\xb_2\xb_3\xb_4\xb_5^{-1}}+2\tr{\xb_1\xb_3\xb_2\xb_4\xb_5^{-1}}.\label{6in5ord}
\end{equation}
In other words, we can always assume there is an ordering on the second two letters, and 
we need only the generator that satisfies this ordering.  This provides us with $12$ generators for every choice for $\xb_5^{-1}$.  We need five further reductions.  However, we now fix a set of twelve to work with
\begin{center}
\begin{align}\label{fixed12}
&\{\T(1, 2, 4, 3, -5), \T(1, 2, 3, 4, -5), \T(1, 3, 4, 2, -5), \T(2, 1, 4, 3, -5),\nonumber\\ 
&\T(2, 1, 3, 4, -5), \T(2, 3, 4, 1, -5), \T(4, 1, 3, 2, -5), \T(4, 2, 3, 1, -5), \nonumber\\
&\T(4, 1, 2, 3, -5), \T(3, 1, 2, 4, -5), \T(3, 1, 4, 2, -5), \T(3, 2, 4, 1, -5)\}. 
\end{align}
\end{center}
Note that if we derive a relation in general it will not be in these terms, and switching to these generators introduces a sign.  For instance, $\T(2, 4, 3, 1, -5)$ must be replaced by $-\T(2, 3, 4, 1, -5)$ since their sum is a reduction relation.

Having found $12$ reductions of the $17$ needed, we must find $5$ more.  Identifying $\wb=\zb$ in the fundamental relation $$\tr{\ub\vb\wb\xb\yb\zb}+\tr{\ub\vb\wb\yb\xb\zb}+\tr{\vb\ub\wb\xb\yb\zb}+\tr{\vb\ub\wb\yb\xb\zb},$$ gives a relation for
$$\tr{\ub\vb\wb\xb\yb\wb}+\tr{\ub\vb\wb\yb\xb\wb}+\tr{\vb\ub\wb\xb\yb\wb}+\tr{\vb\ub\wb\yb\xb\wb},$$ and subsequently using Equation $\eqref{eq:2}$ yields a reduction for
\begin{align*}
&\T(1, 2, 3, 4, -5) + \T(1, 2, 4, 3, -5) + \T(2, 1, 3, 4, -5) + \\
&  \T(2, 1, 4, 3, -5) + \T(3, 4, 1, 2, -5) + \T(3, 4, 2, 1, -5) + \\
&  \T(4, 3, 1, 2, -5) + \T(4, 3, 2, 1, -5).\end{align*}

Imposing the ordering of the second two letters from the reduction Expression $\eqref{6in5ord}$ and changing signs accordingly yields the following expression having trace degree at most $5$ and in terms of our chosen $12$ generators of this type:
\begin{align}\label{eq:r1}
&\T(1, 2, 3, 4, -5) + \T(1, 2, 4, 3, -5) + \T(2, 1, 3, 4, -5) + \nonumber\\
&  \T(2, 1, 4, 3, -5) - \T(3, 1, 4, 2, -5) - \T(3, 2, 4, 1, -5) - \nonumber\\
&  \T(4, 1, 3, 2, -5) - \T(4, 2, 3, 1, -5).\end{align}

Permuting the the letters $\xb_1$ and $\xb_3$ and again putting the result in terms of the $12$ generators of this type provides for a reduction formula for the expression:
\begin{align}\label{eq:r2}
&-\T(1, 2, 4, 3, -5) - \T(1, 3, 4, 2, -5) - \T(2, 1, 3, 4, -5) + \nonumber\\
&  \T(2, 3, 4, 1, -5) - \T(3, 1, 2, 4, -5) + \T(3, 2, 4, 1, -5) + \nonumber\\
&  \T(4, 1, 2, 3, -5) + \T(4, 1, 3, 2, -5).
\end{align}
We will show Expressions $\eqref{eq:r1}$ and $\eqref{eq:r2}$ to be independent in what follows.  But first we need three more reductions relations.

Considering Equation $\eqref{cyclicsum}$, we have a reduction for $$\tr{\xb_2\xb_5\xb_3\sum\xb_5\xb_4\xb_1}$$ which is otherwise expressed
\begin{align*}
&\T(2,5,3,5,4,1)+\T(2,5,3,5,1,4)+\T(2,5,3,4,5,1)\\
&+ \T(2,5,3,4,1,5)+\T(2,5,3,1,5,4)+\T(2,5,3,1,4,5).
\end{align*}
Using Equation $\eqref{eq:2}$, we have a reduction for
\begin{align*}
&-\T(2,5,5,3,4,1)-\T(2,3,5,5,4,1)-\T(2,5,5,3,1,4)-\T(2,3,5,5,1,4)\\ &-\T(2,5,5,3,4,1)-\T(2,3,4,5,5,1)-\T(2,5,5,3,4,1)-\T(2,3,4,1,5,5)\\
&-\T(2,5,5,3,1,4)- \T(2,3,1,5,5,4)-\T(2,5,5,3,1,4)-\T(2,3,1,4,5,5).
\end{align*}

Then switching squares for inverses and cyclically permuting letters so the letter with a negative exponent is at the end of the word gives the expression
\begin{align*}
&-\T(3,4,1,2, -5)- \T(4, 1, 2, 3, -5)- \T(3,1,4,2, -5)\\ 
&- \T(1, 4, 2, 3, -5) -\T(3,4,1,2, -5)-\T(1, 2, 3, 4, -5)\\
&- \T(3,4,1,2, -5)- \T(2, 3, 4, 1, -5)-\T(3,1,4,2, -5)\\
&-\T(4, 2, 3, 1, -5)-\T(3,1,4,2, -5)-\T(2, 3, 1, 4, -5).
\end{align*}

Multiplying by $-1$ and combining like terms then gives a reduction for
\begin{align*}
&\T(1, 2, 3, 4, -5) + \T(1, 4, 2, 3, -5) + \T(2, 3, 1, 4, -5) + \T(2, 3, 4, 1, -5)\\
&+ \T(4, 1, 2, 3, -5) + \T(4, 2, 3, 1, -5)+ 3 \T(3,1,4,2, -5) +3 \T(3,4,1,2, -5).
\end{align*}
However, $\T(3,1,4,2, -5) + \T(3,4,1,2, -5)$ itself has trace degree $5$; that is, we have already shown that this expression is itself entirely reducible (it allowed us to pick an ordering and gave the first $12$ reductions).  Consequently, we have, after re-writting this expression in terms of our chosen $12$ generators of this type and accounting for signs, a reduction for
\begin{align}\label{eq:r3}
&\T(1, 2, 3, 4, -5) - \T(1, 2, 4, 3, -5) - \T(2, 1, 3, 4, -5) +\nonumber\\ 
&  \T(2, 3, 4, 1, -5) + \T(4, 1, 2, 3, -5) + \T(4, 2, 3, 1, -5).
\end{align}

Applying the permutations $(13)$ and $(24)$ to this expression, where the permutation acts on the indices of the generic matrices, and then putting the result in terms of our $12$ chosen generators provides two additional expressions of this type:

\begin{align}\label{eq:r4}
&\T(2, 1, 3, 4, -5) + \T(2, 1, 4, 3, -5) - \T(3, 1, 2, 4, -5) -\nonumber\\ 
&  \T(3, 2, 4, 1, -5) - \T(4, 1, 2, 3, -5) - \T(4, 2, 3, 1, -5),
\end{align}

\begin{align}\label{eq:r5}
&\T(1, 2, 4, 3, -5) - \T(1, 3, 4, 2, -5) + \T(2, 1, 4, 3, -5) -\nonumber\\ 
&  \T(2, 3, 4, 1, -5) - \T(4, 1, 3, 2, -5) - \T(4, 2, 3, 1, -5).
\end{align}

The five relations which correspond to the Expressions $\eqref{eq:r1}, \eqref{eq:r2}, \eqref{eq:r3}, \eqref{eq:r4}, \eqref{eq:r5}$ determine the rows of a $5\times 12$ matrix of $1$'s, $-1$'s, and $0$'s.  Using {\it Mathematica} we compute its rank to be $5$, which implies the relations are independent and we may make five further reductions of the $12$ generators, leaving only $7$, as required.  Again, that is $7$ choices for every choice for the last letter in the word to have an exponent of $-1$; there are $5$ such choices.  Hence we have $35$ generators of this type; the minimal number available.

With respect to the ordering in the list of $12$ generators in $\eqref{fixed12}$, the $5\times 12$ matrix is:

{\tiny
$$
\left(
\begin{array}{llllllllllll}
 1 & 1 & 0 & 1 & 1 & 0 & -1 & -1 & 0 & 0 & -1 & -1\\
-1 & 0 & -1 & 0 & -1 & 1 & 1 & 0 & 1 & -1 & 0 & 1\\
 -1 & 1 & 0 & 0 & -1 & 1 & 0 & 1 & 1 & 0 & 0 & 0\\
 0 & 0 & 0 & 1 & 1 & 0 & 0 & -1 & -1 & -1 & 0 & -1\\
 1 & 0 & -1 & 1 & 0 & -1 & -1 & -1 & 0 & 0 & 0 & 0\\
\end{array}
\right)$$}

One such solution set of seven among the twelve generators with $\xb_5^{-1}$ as its last letter is 
\begin{align*}
&\{\T(2, 3, 4, 1, -5), \T(3, 1, 2, 4, -5), \T(3, 1, 4, 2, -5),\\
&\T(3, 2, 4, 1, -5),\T(4, 1, 2, 3, -5), \T(4, 1, 3, 2, -5), \T(4, 2, 3, 1, -5)\}.
\end{align*}  
For the other $28$ apply each of the permutations $(15),(25),(35),$ and $(45)$ to this set, one at a time, to get four more sets of $7$.

The last case to consider is the case of six independent generic matrices and generators of length $6$ in those terms.  The only generator type not accounted for yet is $\tr{\ub\vb\wb\xb\yb\zb}$.

The table is as follows:

\begin{table}[!ht]
\begin{center}
\begin{tabular}{|c|c|}
\hline
Minimal Number & Generators\\
\hline
$\N{6,6}=15\binom{r}{6}$ & $\tr{\ub\vb\wb\xb\yb\zb}$\\
\hline
\end{tabular}
\end{center}
\caption{This table lists the minimal generators in words with only six letters.}\label{min6}
\end{table}
The weights for the irreps corresponding to the generators of length $6$ in $6$ letters are, like in the last case, $\{(3,1,1,1),(3,3)\}$; but in this case the total of their dimensions is $770$.  The previous cases $r=1,2,3,4,5$ account for $35\binom{6}{5}+26\binom{6}{4} + 7\binom{6}{3} + \binom{6}{2}=755$ of these generators.  This leaves $15$ of type $\tr{\ub\vb\wb\xb\yb\zb}$.  Cyclically there are $120$ possibilities.  We must provide $105$ independent relations.

In \cite{MS1} similar relations for $\C[\Y_6]$ were given. They were of two types: $37$ of $\tr{\ub\sum\xb\yb\zb}$ and $68$ of  $$\tr{\ub\vb\wb\xb\yb\zb}+\tr{\ub\vb\wb\yb\xb\zb}+\tr{\vb\ub\wb\xb\yb\zb}+\tr{\vb\ub\wb\yb\xb\zb}.$$

However, we construct the $105$ reductions from two variations of the relation for $\tr{\ub\sum\xb\yb\zb}$ alone.

First we have
\begin{align}\label{6in6eq:1}
& \T(1,2,3,4,5,6)+\T(1,2,3,5,4,6)+\T(4,1,2,3,5,6)\nonumber\\
&+\T(4,5,1,2,3,6)+\T(5,1,2,3,4,6)+\T(5,4,1,2,3,6)=\nonumber\\
&\nonumber\\
&\T(4) \T(5) \T(6) \T(1,2,3)-\T(6) \T(4,5) \T(1,2,3)-\T(5) \T(4,6) \T(1,2,3)\nonumber\\
&-\T(4) \T(5,6) \T(1,2,3)+\T(5,4,6)\T(1,2,3)+\T(5,6,4) \T(1,2,3)\nonumber\\
&-\T(5) \T(6) \T(1,2,3,4)+\T(5,6) \T(1,2,3,4)-\T(4) \T(6) \T(1,2,3,5)\nonumber\\
&+\T(4,6)\T(1,2,3,5)-\T(4) \T(5) \T(1,2,3,6)+\T(4,5) \T(1,2,3,6)\nonumber\\
&+\T(6) \T(1,2,3,4,5)+\T(5) \T(1,2,3,4,6)+\T(6)\T(1,2,3,5,4)\nonumber\\
&+\T(4) \T(1,2,3,5,6)+\T(5) \T(1,2,3,6,4)+\T(4) \T(1,2,3,6,5),
\end{align}
from setting $\ub=\xb_1\xb_2\xb_3,\xb=\xb_4,\yb=\xb_5, \zb=\xb_6$ in $\tr{\ub\sum\xb\yb\zb}$.

In what follows we will cyclically permute all generators of this type (there are 720 of them) so that the generic matrix $\xb_6$ is at the end of the word.  Then we have $120$ generators in the following natural order:
{\small
\begin{align}\label{6in6gens}
&\{\T(1,2,3,4,5,6),\T(1,2,3,5,4,6),\T(1,2,4,3,5,6),\T(1,2,4,5,3,6),\T(1,2,5,3,4,6),\nonumber\\
&\T(1,2,5,4,3,6),\T(1,3,2,4,5,6),\T(1,3,2,5,4,6),\T(1,3,4,2,5,6),\T(1,3,4,5,2,6),\nonumber\\
&\T(1,3,5,2,4,6),\T(1,3,5,4,2,6),\T(1,4,2,3,5,6),\T(1,4,2,5,3,6),\T(1,4,3,2,5,6),\nonumber\\
&\T(1,4,3,5,2,6),\T(1,4,5,2,3,6),\T(1,4,5,3,2,6),\T(1,5,2,3,4,6),\T(1,5,2,4,3,6),\nonumber\\
&\T(1,5,3,2,4,6),\T(1,5,3,4,2,6),\T(1,5,4,2,3,6),\T(1,5,4,3,2,6),\T(2,1,3,4,5,6),\nonumber\\
&\T(2,1,3,5,4,6),\T(2,1,4,3,5,6),\T(2,1,4,5,3,6),\T(2,1,5,3,4,6),\T(2,1,5,4,3,6),\nonumber\\
&\T(2,3,1,4,5,6),\T(2,3,1,5,4,6),\T(2,3,4,1,5,6),\T(2,3,4,5,1,6),\T(2,3,5,1,4,6),\nonumber\\
&\T(2,3,5,4,1,6),\T(2,4,1,3,5,6),\T(2,4,1,5,3,6),\T(2,4,3,1,5,6),\T(2,4,3,5,1,6),\nonumber\\
&\T(2,4,5,1,3,6),\T(2,4,5,3,1,6),\T(2,5,1,3,4,6),\T(2,5,1,4,3,6),\T(2,5,3,1,4,6),\nonumber\\
&\T(2,5,3,4,1,6),\T(2,5,4,1,3,6),\T(2,5,4,3,1,6),\T(3,1,2,4,5,6),\T(3,1,2,5,4,6),\nonumber\\
&\T(3,1,4,2,5,6),\T(3,1,4,5,2,6),\T(3,1,5,2,4,6),\T(3,1,5,4,2,6),\T(3,2,1,4,5,6),\nonumber\\
&\T(3,2,1,5,4,6),\T(3,2,4,1,5,6),\T(3,2,4,5,1,6),\T(3,2,5,1,4,6),\T(3,2,5,4,1,6),\nonumber\\
&\T(3,4,1,2,5,6),\T(3,4,1,5,2,6),\T(3,4,2,1,5,6),\T(3,4,2,5,1,6),\T(3,4,5,1,2,6),\nonumber\\
&\T(3,4,5,2,1,6),\T(3,5,1,2,4,6),\T(3,5,1,4,2,6),\T(3,5,2,1,4,6),\T(3,5,2,4,1,6),\nonumber\\
&\T(3,5,4,1,2,6),\T(3,5,4,2,1,6),\T(4,1,2,3,5,6),\T(4,1,2,5,3,6),\T(4,1,3,2,5,6),\nonumber\\
&\T(4,1,3,5,2,6),\T(4,1,5,2,3,6),\T(4,1,5,3,2,6),\T(4,2,1,3,5,6),\T(4,2,1,5,3,6),\nonumber\\
&\T(4,2,3,1,5,6),\T(4,2,3,5,1,6),\T(4,2,5,1,3,6),\T(4,2,5,3,1,6),\T(4,3,1,2,5,6),\nonumber\\
&\T(4,3,1,5,2,6),\T(4,3,2,1,5,6),\T(4,3,2,5,1,6),\T(4,3,5,1,2,6),\T(4,3,5,2,1,6),\nonumber\\
&\T(4,5,1,2,3,6),\T(4,5,1,3,2,6),\T(4,5,2,1,3,6),\T(4,5,2,3,1,6),\T(4,5,3,1,2,6),\nonumber\\
&\T(4,5,3,2,1,6),\T(5,1,2,3,4,6),\T(5,1,2,4,3,6),\T(5,1,3,2,4,6),\T(5,1,3,4,2,6),\nonumber\\
&\T(5,1,4,2,3,6),\T(5,1,4,3,2,6),\T(5,2,1,3,4,6),\T(5,2,1,4,3,6),\T(5,2,3,1,4,6),\nonumber\\
&\T(5,2,3,4,1,6),\T(5,2,4,1,3,6),\T(5,2,4,3,1,6),\T(5,3,1,2,4,6),\T(5,3,1,4,2,6),\nonumber\\
&\T(5,3,2,1,4,6),\T(5,3,2,4,1,6),\T(5,3,4,1,2,6),\T(5,3,4,2,1,6),\T(5,4,1,2,3,6),\nonumber\\
&\T(5,4,1,3,2,6),\T(5,4,2,1,3,6),\T(5,4,2,3,1,6),\T(5,4,3,1,2,6),\T(5,4,3,2,1,6)\}.
\end{align}} 
Now by applying the $120$ permutations of the five letters $\xb_1,...,\xb_5$ to Equation $\eqref{6in6eq:1}$, we come to exactly 60 new relations.  Note that this permutation action amounts to just relabeling the letters since these relations are determined by words of the given form; the labels are insignificant by themselves.

For the second variation of $\tr{\ub\sum\xb\yb\zb}$, let $\ub=\xb_1,\xb=\xb_2\xb_3,\yb=\xb_4\xb_5, \zb=\xb_6$.  Then we come to the relation

\begin{align}\label{6in6eq:2}
 &\T(1,2,3,4,5,6)+\T(1,4,5,2,3,6)+\T(2,3,1,4,5,6)\nonumber\\
&+\T(2,3,4,5,1,6)+\T(4,5,1,2,3,6)+\T(4,5,2,3,1,6)=\nonumber\\
&\nonumber\\
&\T(1) \T(6) \T(2,3) \T(4,5)-\T(1,6) \T(2,3) \T(4,5)-\T(6) \T(1,2,3) \T(4,5)\nonumber\\
&-\frac{1}{2} \T(1) \T(2,3,6)\T(4,5)-\frac{1}{2} \T(1) \T(6,2,3) \T(4,5)
+\T(1,2,3,6) \T(4,5)\nonumber\\
&+\T(1,6,2,3) \T(4,5)-\T(6) \T(2,3)\T(1,4,5)+\frac{1}{2} \T(1,4,5) \T(2,3,6)\nonumber\\
&+\frac{1}{2} \T(1,4,5) \T(6,2,3)-\T(1) \T(2,3) \T(6,4,5)+\T(1,2,3)\T(6,4,5)\nonumber\\
&+\T(2,3) \T(1,4,5,6)+\T(2,3) \T(1,6,4,5)-\T(1) \T(6) \T(2,3,4,5)\nonumber\\
&+\T(1,6) \T(2,3,4,5)+\T(6)\T(1,2,3,4,5)+\T(6) \T(1,4,5,2,3)\nonumber\\
&+\T(1) \T(4,5,2,3,6)+\T(1) \T(4,5,6,2,3).
\end{align}
Again applying the $120$ permutations of the first five generic matrices to Equation $\eqref{6in6eq:2}$ yields exactly $60$ equations; for a total of $120$ equations when added to the $60$ equations coming from the permutations of Equation $\eqref{6in6eq:1}$.

The left-hand-sides of these relations are linear in the $120$ generators of length $6$ with coefficients only $1$ or $0$.

With respect to the order of the generators implicit in the list $\eqref{6in6gens}$ above, we form the $120\times 120$ matrix:
$$M=\left(
\begin{array}{lll}
 M_{11} & M_{12} & M_{13} \\
M_{21} & M_{22} & M_{23} \\
M_{31} & M_{32} & M_{33} \\
\end{array}\right),$$ 
where each entry is a $40\times 40$ matrix of ones and zeros.  The rows of $M$ correspond to the $120$ relations coming from the permutations of Equations $\eqref{6in6eq:1}$ and $\eqref{6in6eq:2}$.  We list the submatrices $M_{ij}$ in the appendix for completeness.

Computing the rank of $M$ with {\it Mathematica} we determine it to be $105$.  Hence we can reduce the $120$ generators to exactly $15$; the known minimal number.

The complement of any $105$ pivot columns will correspond to $15$ minimal generators.  For instance, removing columns $86, 88,90,104,108,110-114$, and $116-120$ gives a full rank matrix.  Hence the generators which correspond to those columns (the $86^{\text{th}}, 88^{\text{th}}, ...,$ etc. entries in the above list $\eqref{6in6gens}$ of $120$ generators) are a minimal set of generators of this type.  Here they are:

\begin{align*}
& \{\T(4,3,1,5,2,6),\T(4,3,2,5,1,6),\T(4,3,5,2,1,6),\T(5,2,1,4,3,6),\\
&\T(5,2,4,3,1,6),\T(5,3,1,4,2,6),\T(5,3,2,1,4,6),\T(5,3,2,4,1,6),\\
&\T(5,3,4,1,2,6),\T(5,3,4,2,1,6),\T(5,4,1,3,2,6),\T(5,4,2,1,3,6),\\
&\T(5,4,2,3,1,6),\T(5,4,3,1,2,6),\T(5,4,3,2,1,6)\}.
\end{align*}

Performing row reductions on the $120\times 120$ matrix and using the permutations and explicit relations above, one can find the exact relations necessary to remove these $105$ generators (this is not computationally trivial however).

This concludes the proof of Theorem $\ref{maintheorem}$.

\section*{Acknowledgments}
The author thanks Bill Goldman for introducing him to this topic and for his support and encouragement. This work was completed at Kansas State University during the 2006-2007 academic year, but the writing was completed at Instituto Superior T\'ecnico during the 2007-2008 academic year.  The author thanks both KSU and IST for providing stimulating environments to work.  We thank Artem Lopatin for helpful references and comments, and we thank Mara Neusel, Frank Grosshans, Zongshu Lin, Gustavo Granja, Pierre Will, and Adam Sikora for general interest in this work.  Also we thank an anonymous referee for helpful comments and suggestions.

\section*{Appendix A}

The algorithm alluded to at the end of the proof of Lemma \ref{generatorform} is as follows (see \cite{L2} for details):

STEP $1$:  Define $\pol_2(\xb,\yb,\zb)=$ $$\pol(\yb,\xb^2\zb)+\xb\pol(\yb,\xb\zb)-\pol(\xb,\yb^2)\zb-\pol(\xb,\yb)\zb\yb+\xb^2\pol(\yb,\zb).$$  Then $3\xb^2\zb\yb^2=\pol_2(\xb,\yb,\zb)$.

STEP $2$:  Define $\mathrm{pre}\pol_3(\xb,\ub,\vb,\zb)=$ $$\pol_2(\xb,\ub+\vb,\zb)-\pol_2(\xb,\ub,\zb)-\pol_2(\xb,\vb,\zb),$$ and define 
$\pol_3(\xb,\ub,\vb,\wb,\zb)=$ $$\mathrm{pre}\pol_3(\xb,\ub,\vb,\zb\wb)+\mathrm{pre}\pol_3(\xb,\wb\ub,\vb,\zb)-\mathrm{pre}\pol_3(\xb,\wb,\vb,\zb)\ub.$$  
Then $6\xb^2\zb\wb\ub\vb=\pol_3(\xb,\ub,\vb,\wb,\zb)$.

STEP $3$:  $6(\xb\yb\zb\wb\ub\vb+\yb\xb\zb\wb\ub\vb)=$ $$\pol_3(\xb+\yb,\ub,\vb,\wb,\zb)-\pol_3(\xb,\ub,\vb,\wb,\zb)-\pol_3(\yb,\ub,\vb,\wb,\zb).$$
NOTE: In \cite{L2} the expression on the right in STEP $3$ is shown to have trace degree $6$ and total degree $5$ (the total degree is the largest word in the expression that is not the argument of a trace).

STEP $4$:  Using Mathematica, we implement STEPS $1$-$3$ to find:
\begin{align*}&6(\xb\yb\zb\wb\ub\vb+\yb\xb\zb\wb\ub\vb)=\\\\
&\id\bigg(2\tr{\xb\yb\zb\wb\ub\vb+\yb\xb\zb\wb\ub\vb}\\
&+\tr{\xb\yb\zb\wb\vb\ub+\xb\yb\zb\vb\wb\ub+\yb\xb\zb\wb\vb\ub+\yb\xb\zb\vb\wb\ub}\bigg)\\
&+\mathrm{E}(\xb,\yb,\zb,\wb,\ub,\vb),
\end{align*}
where $\mathrm{E}$ is a polynomial expression of total degree and trace degree less than or equal to $5$.

STEP $5$:  Take the trace of both sides of the expression in STEP $4$ and cancel like terms to yield: 
$$3\bigg(\tr{\xb\yb\zb\wb\vb\ub}+\tr{\xb\yb\zb\vb\wb\ub}+\tr{\yb\xb\zb\wb\vb\ub}+\tr{\yb\xb\zb\vb\wb\ub}\bigg)=$$ 
$\mathrm{F}(\xb,\yb,\zb,\wb,\vb,\ub)$, where $\mathrm{F}$ is a polynomial expression of trace degree less than or equal to $5$.

This gives the desired expression which we refer to as the {\it fundamental relation}.

\section*{Appendix B}
With respect to the order of the generators implicit in the list $\eqref{6in6gens}$ we form the $120\times 120$ matrix:
$$M=\left(
\begin{array}{lll}
 M_{11} & M_{12} & M_{13} \\
M_{21} & M_{22} & M_{23} \\
M_{31} & M_{32} & M_{33} \\
\end{array}\right),$$ 
where each entry is a $40\times 40$ matrix of ones and zeros.  The rows of $M$ correspond to the $120$ relations coming from the permutations of Equations $\eqref{6in6eq:1}$ and $\eqref{6in6eq:2}$.  We list the submatrices $M_{ij}$ here for completeness:
\newpage
$M_{11}=$\\

\hspace{-4.7cm}
{\tiny
$\left(
\begin{array}{llllllllllllllllllllllllllllllllllllllll}
 0 & 0 & 0 & 0 & 0 & 0 & 0 & 0 & 0 & 0 & 0 & 0 & 0 & 0 & 0 & 0 & 1 & 1 & 0 & 0 & 0 & 0 & 0 & 0 & 0 & 0
   & 0 & 1 & 0 & 0 & 1 & 0 & 0 & 0 & 0 & 0 & 0 & 0 & 0 & 0 \\
 0 & 0 & 0 & 0 & 0 & 0 & 0 & 0 & 0 & 0 & 0 & 0 & 0 & 0 & 0 & 0 & 0 & 0 & 0 & 0 & 0 & 0 & 1 & 1 & 0 & 0
   & 0 & 0 & 0 & 1 & 0 & 1 & 0 & 0 & 0 & 0 & 0 & 0 & 0 & 0 \\
 0 & 0 & 0 & 0 & 0 & 0 & 0 & 0 & 0 & 0 & 0 & 0 & 0 & 0 & 0 & 0 & 0 & 0 & 0 & 1 & 1 & 0 & 0 & 0 & 0 & 0
   & 0 & 0 & 0 & 0 & 0 & 0 & 0 & 0 & 0 & 0 & 0 & 1 & 1 & 0 \\
 0 & 0 & 0 & 1 & 0 & 0 & 1 & 0 & 0 & 0 & 0 & 0 & 0 & 0 & 0 & 0 & 0 & 0 & 0 & 0 & 0 & 0 & 0 & 0 & 0 & 0
   & 0 & 0 & 0 & 0 & 0 & 0 & 0 & 0 & 0 & 0 & 0 & 0 & 0 & 0 \\
 0 & 0 & 0 & 0 & 0 & 0 & 0 & 0 & 0 & 0 & 0 & 0 & 0 & 1 & 1 & 0 & 0 & 0 & 0 & 0 & 0 & 0 & 0 & 0 & 0 & 0
   & 0 & 0 & 0 & 0 & 0 & 0 & 0 & 0 & 0 & 0 & 0 & 0 & 0 & 0 \\
 0 & 0 & 0 & 0 & 0 & 1 & 0 & 1 & 0 & 0 & 0 & 0 & 0 & 0 & 0 & 0 & 0 & 0 & 0 & 0 & 0 & 0 & 0 & 0 & 0 & 0
   & 0 & 0 & 0 & 0 & 0 & 0 & 0 & 0 & 0 & 0 & 0 & 0 & 0 & 0 \\
 0 & 0 & 0 & 0 & 0 & 0 & 0 & 0 & 0 & 0 & 0 & 0 & 0 & 0 & 0 & 0 & 0 & 0 & 1 & 0 & 0 & 1 & 0 & 0 & 0 & 0
   & 0 & 0 & 1 & 0 & 0 & 0 & 1 & 0 & 0 & 0 & 0 & 0 & 0 & 0 \\
 0 & 0 & 0 & 0 & 1 & 0 & 0 & 0 & 1 & 0 & 0 & 0 & 0 & 0 & 0 & 0 & 0 & 0 & 0 & 0 & 0 & 0 & 0 & 0 & 0 & 0
   & 0 & 0 & 0 & 0 & 0 & 0 & 0 & 0 & 0 & 0 & 0 & 0 & 0 & 0 \\
\end{array}
\right)$}

\newpage
$M_{12}=$\\

\hspace{-4.7cm}
{\tiny $\left(
\begin{array}{llllllllllllllllllllllllllllllllllllllll}
 0 & 0 & 0 & 0 & 0 & 0 & 0 & 0 & 0 & 0 & 0 & 1 & 0 & 0 & 1 & 0 & 0 & 0 & 0 & 0 & 0 & 0 & 0 & 0 & 0 & 0
   & 0 & 0 & 0 & 0 & 0 & 0 & 0 & 0 & 0 & 0 & 0 & 0 & 0 & 0 \\
 0 & 0 & 0 & 0 & 0 & 0 & 0 & 0 & 0 & 0 & 0 & 0 & 0 & 1 & 0 & 1 & 0 & 0 & 0 & 0 & 0 & 0 & 0 & 0 & 0 & 0
   & 0 & 0 & 0 & 0 & 0 & 0 & 0 & 0 & 0 & 0 & 0 & 0 & 0 & 0 \\
 0 & 0 & 0 & 0 & 0 & 0 & 0 & 0 & 0 & 0 & 0 & 0 & 1 & 0 & 0 & 0 & 1 & 0 & 0 & 0 & 0 & 0 & 0 & 0 & 0 & 0
   & 0 & 0 & 0 & 0 & 0 & 0 & 0 & 0 & 0 & 0 & 0 & 0 & 0 & 0 \\
 1 & 1 & 0 & 0 & 0 & 0 & 0 & 0 & 1 & 0 & 0 & 0 & 0 & 0 & 0 & 0 & 0 & 1 & 0 & 0 & 0 & 0 & 0 & 0 & 0 & 0
   & 0 & 0 & 0 & 0 & 0 & 0 & 0 & 0 & 0 & 0 & 0 & 0 & 0 & 0 \\
 0 & 0 & 0 & 1 & 1 & 0 & 0 & 0 & 0 & 0 & 1 & 0 & 0 & 0 & 0 & 0 & 0 & 0 & 1 & 0 & 0 & 0 & 0 & 0 & 0 & 0
   & 0 & 0 & 0 & 0 & 0 & 0 & 0 & 0 & 0 & 0 & 0 & 0 & 0 & 0 \\
 0 & 0 & 0 & 0 & 0 & 0 & 1 & 1 & 0 & 1 & 0 & 0 & 0 & 0 & 0 & 0 & 0 & 0 & 0 & 1 & 0 & 0 & 0 & 0 & 0 & 0
   & 0 & 0 & 0 & 0 & 0 & 0 & 0 & 0 & 0 & 0 & 0 & 0 & 0 & 0 \\
 0 & 0 & 0 & 0 & 0 & 0 & 0 & 0 & 0 & 0 & 0 & 0 & 0 & 0 & 0 & 0 & 0 & 0 & 0 & 0 & 0 & 1 & 1 & 0 & 0 & 0
   & 0 & 0 & 0 & 0 & 0 & 0 & 0 & 0 & 0 & 0 & 0 & 0 & 0 & 0 \\
 0 & 0 & 1 & 0 & 0 & 1 & 0 & 0 & 0 & 0 & 0 & 0 & 0 & 0 & 0 & 0 & 0 & 0 & 0 & 0 & 1 & 0 & 0 & 1 & 0 & 0
   & 0 & 0 & 0 & 0 & 0 & 0 & 0 & 0 & 0 & 0 & 0 & 0 & 0 & 0 \\
\end{array}
\right)$}

\newpage
$M_{13}=$\\

\hspace{-4.7cm}
{\tiny $ 
\left(
\begin{array}{llllllllllllllllllllllllllllllllllllllll}
 0 & 0 & 0 & 0 & 0 & 0 & 0 & 0 & 0 & 0 & 0 & 0 & 0 & 0 & 0 & 0 & 0 & 0 & 0 & 0 & 0 & 0 & 0 & 0 & 0 & 0
   & 0 & 0 & 0 & 0 & 0 & 0 & 0 & 0 & 0 & 0 & 0 & 0 & 0 & 0 \\
 0 & 0 & 0 & 0 & 0 & 0 & 0 & 0 & 0 & 0 & 0 & 0 & 0 & 0 & 0 & 0 & 0 & 0 & 0 & 0 & 0 & 0 & 0 & 0 & 0 & 0
   & 0 & 0 & 0 & 0 & 0 & 0 & 0 & 0 & 0 & 0 & 0 & 0 & 0 & 0 \\
 0 & 0 & 0 & 0 & 0 & 0 & 0 & 0 & 0 & 0 & 0 & 0 & 0 & 0 & 0 & 0 & 0 & 0 & 0 & 0 & 0 & 0 & 0 & 0 & 0 & 0
   & 0 & 0 & 0 & 0 & 0 & 0 & 0 & 0 & 0 & 0 & 0 & 0 & 0 & 0 \\
 0 & 0 & 0 & 0 & 0 & 0 & 0 & 0 & 0 & 0 & 0 & 0 & 0 & 0 & 0 & 0 & 0 & 0 & 0 & 0 & 0 & 0 & 0 & 0 & 0 & 0
   & 0 & 0 & 0 & 0 & 0 & 0 & 0 & 0 & 0 & 0 & 0 & 0 & 0 & 0 \\
 0 & 0 & 0 & 0 & 0 & 0 & 0 & 0 & 0 & 0 & 0 & 0 & 0 & 0 & 0 & 0 & 0 & 0 & 0 & 0 & 0 & 0 & 0 & 0 & 0 & 0
   & 0 & 0 & 0 & 0 & 0 & 0 & 0 & 0 & 0 & 0 & 0 & 0 & 0 & 0 \\
 0 & 0 & 0 & 0 & 0 & 0 & 0 & 0 & 0 & 0 & 0 & 0 & 0 & 0 & 0 & 0 & 0 & 0 & 0 & 0 & 0 & 0 & 0 & 0 & 0 & 0
   & 0 & 0 & 0 & 0 & 0 & 0 & 0 & 0 & 0 & 0 & 0 & 0 & 0 & 0 \\
 0 & 0 & 0 & 0 & 0 & 0 & 0 & 0 & 0 & 0 & 0 & 0 & 0 & 0 & 0 & 0 & 0 & 0 & 0 & 0 & 0 & 0 & 0 & 0 & 0 & 0
   & 0 & 0 & 0 & 0 & 0 & 0 & 0 & 0 & 0 & 0 & 0 & 0 & 0 & 0 \\
 0 & 0 & 0 & 0 & 0 & 0 & 0 & 0 & 0 & 0 & 0 & 0 & 0 & 0 & 0 & 0 & 0 & 0 & 0 & 0 & 0 & 0 & 0 & 0 & 0 & 0
   & 0 & 0 & 0 & 0 & 0 & 0 & 0 & 0 & 0 & 0 & 0 & 0 & 0 & 0 \\
\end{array}
\right)
$}
\newpage
$ M_{21}=$\\

\hspace{-4.7cm}
{\tiny $ 
\left(
\begin{array}{llllllllllllllllllllllllllllllllllllllll}
 0 & 0 & 0 & 0 & 0 & 0 & 1 & 0 & 0 & 1 & 0 & 0 & 0 & 0 & 0 & 0 & 0 & 0 & 0 & 0 & 0 & 0 & 0 & 0 & 1 & 0
   & 0 & 0 & 0 & 0 & 0 & 0 & 0 & 0 & 0 & 0 & 0 & 0 & 0 & 0 \\
 1 & 0 & 0 & 0 & 0 & 0 & 0 & 0 & 0 & 0 & 0 & 0 & 0 & 0 & 0 & 0 & 1 & 0 & 0 & 0 & 0 & 0 & 0 & 0 & 0 & 0
   & 0 & 0 & 0 & 0 & 1 & 0 & 0 & 1 & 0 & 0 & 0 & 0 & 0 & 0 \\
 0 & 0 & 0 & 0 & 0 & 0 & 0 & 0 & 0 & 1 & 0 & 0 & 0 & 0 & 0 & 0 & 1 & 0 & 0 & 0 & 0 & 0 & 0 & 0 & 0 & 0
   & 0 & 0 & 0 & 0 & 0 & 0 & 0 & 0 & 0 & 0 & 0 & 0 & 0 & 0 \\
 1 & 0 & 0 & 1 & 0 & 0 & 0 & 0 & 0 & 0 & 0 & 0 & 0 & 0 & 0 & 0 & 0 & 0 & 0 & 0 & 0 & 0 & 0 & 0 & 0 & 0
   & 0 & 0 & 0 & 0 & 0 & 0 & 0 & 0 & 0 & 0 & 0 & 0 & 0 & 0 \\
 0 & 0 & 0 & 0 & 0 & 0 & 0 & 0 & 0 & 0 & 0 & 0 & 0 & 0 & 0 & 0 & 0 & 0 & 0 & 0 & 0 & 0 & 0 & 0 & 0 & 0
   & 0 & 0 & 0 & 0 & 1 & 0 & 0 & 0 & 0 & 0 & 0 & 0 & 0 & 0 \\
 0 & 0 & 0 & 0 & 0 & 0 & 1 & 0 & 0 & 0 & 0 & 0 & 0 & 0 & 0 & 0 & 0 & 1 & 0 & 0 & 0 & 0 & 0 & 0 & 0 & 0
   & 0 & 0 & 0 & 0 & 0 & 0 & 0 & 0 & 0 & 0 & 0 & 0 & 0 & 0 \\
 0 & 0 & 0 & 0 & 0 & 0 & 0 & 0 & 0 & 0 & 0 & 0 & 0 & 0 & 0 & 0 & 0 & 0 & 0 & 0 & 0 & 0 & 0 & 0 & 1 & 0
   & 0 & 1 & 0 & 0 & 0 & 0 & 0 & 0 & 0 & 0 & 0 & 0 & 0 & 0 \\
 0 & 0 & 0 & 1 & 0 & 0 & 0 & 0 & 0 & 0 & 0 & 0 & 0 & 0 & 0 & 0 & 0 & 1 & 0 & 0 & 0 & 0 & 0 & 0 & 0 & 0
   & 0 & 1 & 0 & 0 & 0 & 0 & 0 & 0 & 0 & 0 & 0 & 0 & 0 & 0 \\
\end{array}
\right)
$}
\newpage

$M_{22}=$\\

\hspace{-4.7cm}
{\tiny $
\left(
\begin{array}{llllllllllllllllllllllllllllllllllllllll}
 1 & 0 & 0 & 0 & 0 & 0 & 0 & 0 & 0 & 0 & 0 & 0 & 0 & 0 & 0 & 0 & 0 & 0 & 0 & 0 & 0 & 0 & 0 & 0 & 0 & 0
   & 0 & 0 & 0 & 0 & 0 & 0 & 0 & 0 & 0 & 0 & 0 & 0 & 0 & 0 \\
 0 & 0 & 0 & 0 & 0 & 0 & 0 & 0 & 0 & 0 & 0 & 0 & 0 & 0 & 0 & 0 & 0 & 0 & 0 & 0 & 0 & 0 & 0 & 0 & 0 & 0
   & 0 & 0 & 0 & 0 & 0 & 0 & 0 & 0 & 0 & 0 & 0 & 0 & 0 & 0 \\
 0 & 0 & 0 & 0 & 0 & 0 & 0 & 0 & 0 & 0 & 0 & 1 & 0 & 0 & 0 & 0 & 0 & 0 & 0 & 0 & 0 & 0 & 0 & 0 & 0 & 1
   & 0 & 0 & 0 & 0 & 0 & 0 & 0 & 0 & 0 & 0 & 0 & 0 & 0 & 0 \\
 0 & 0 & 0 & 0 & 0 & 0 & 0 & 0 & 1 & 0 & 0 & 0 & 0 & 0 & 0 & 0 & 0 & 0 & 0 & 0 & 0 & 0 & 0 & 0 & 1 & 0
   & 0 & 0 & 0 & 0 & 0 & 0 & 0 & 0 & 0 & 0 & 0 & 0 & 0 & 0 \\
 0 & 1 & 0 & 0 & 0 & 0 & 0 & 0 & 1 & 0 & 0 & 1 & 0 & 0 & 0 & 0 & 0 & 0 & 0 & 0 & 0 & 0 & 0 & 0 & 0 & 0
   & 0 & 0 & 0 & 0 & 0 & 0 & 0 & 0 & 0 & 0 & 0 & 0 & 0 & 0 \\
 0 & 0 & 0 & 0 & 0 & 0 & 0 & 0 & 0 & 0 & 0 & 0 & 0 & 0 & 1 & 0 & 0 & 1 & 0 & 0 & 0 & 0 & 0 & 0 & 0 & 0
   & 0 & 0 & 0 & 0 & 0 & 0 & 0 & 0 & 0 & 0 & 0 & 0 & 0 & 0 \\
 0 & 0 & 0 & 0 & 0 & 0 & 0 & 0 & 0 & 0 & 0 & 0 & 0 & 0 & 1 & 0 & 0 & 0 & 0 & 0 & 0 & 0 & 0 & 0 & 0 & 1
   & 0 & 0 & 0 & 0 & 0 & 0 & 0 & 0 & 0 & 0 & 0 & 0 & 0 & 0 \\
 0 & 1 & 0 & 0 & 0 & 0 & 0 & 0 & 0 & 0 & 0 & 0 & 0 & 0 & 0 & 0 & 0 & 0 & 0 & 0 & 0 & 0 & 0 & 0 & 0 & 0
   & 0 & 0 & 0 & 0 & 0 & 0 & 0 & 0 & 0 & 0 & 0 & 0 & 0 & 0 \\
\end{array}
\right)$}
\newpage

$M_{23}=$\\

\hspace{-4.7cm}
{\tiny $
\left(
\begin{array}{llllllllllllllllllllllllllllllllllllllll}
 0 & 0 & 0 & 0 & 0 & 0 & 0 & 0 & 0 & 0 & 0 & 1 & 1 & 0 & 0 & 0 & 0 & 0 & 0 & 0 & 0 & 0 & 0 & 0 & 0 & 0
   & 0 & 0 & 0 & 0 & 0 & 0 & 0 & 0 & 0 & 0 & 0 & 0 & 0 & 0 \\
 0 & 0 & 0 & 0 & 0 & 0 & 0 & 0 & 0 & 0 & 1 & 0 & 0 & 1 & 0 & 0 & 0 & 0 & 0 & 0 & 0 & 0 & 0 & 0 & 0 & 0
   & 0 & 0 & 0 & 0 & 0 & 0 & 0 & 0 & 0 & 0 & 0 & 0 & 0 & 0 \\
 0 & 0 & 0 & 0 & 0 & 0 & 0 & 0 & 0 & 0 & 0 & 0 & 1 & 1 & 0 & 0 & 0 & 0 & 0 & 0 & 0 & 0 & 0 & 0 & 0 & 0
   & 0 & 0 & 0 & 0 & 0 & 0 & 0 & 0 & 0 & 0 & 0 & 0 & 0 & 0 \\
 0 & 0 & 0 & 0 & 0 & 0 & 0 & 0 & 0 & 0 & 1 & 0 & 0 & 0 & 1 & 0 & 0 & 0 & 0 & 0 & 0 & 0 & 0 & 0 & 0 & 0
   & 0 & 0 & 0 & 0 & 0 & 0 & 0 & 0 & 0 & 0 & 0 & 0 & 0 & 0 \\
 0 & 0 & 0 & 0 & 0 & 0 & 0 & 0 & 0 & 0 & 0 & 0 & 0 & 1 & 1 & 0 & 0 & 0 & 0 & 0 & 0 & 0 & 0 & 0 & 0 & 0
   & 0 & 0 & 0 & 0 & 0 & 0 & 0 & 0 & 0 & 0 & 0 & 0 & 0 & 0 \\
 0 & 0 & 0 & 0 & 0 & 0 & 0 & 0 & 0 & 0 & 0 & 1 & 0 & 0 & 0 & 1 & 0 & 0 & 0 & 0 & 0 & 0 & 0 & 0 & 0 & 0
   & 0 & 0 & 0 & 0 & 0 & 0 & 0 & 0 & 0 & 0 & 0 & 0 & 0 & 0 \\
 0 & 0 & 0 & 0 & 0 & 0 & 0 & 0 & 0 & 0 & 0 & 0 & 1 & 0 & 0 & 1 & 0 & 0 & 0 & 0 & 0 & 0 & 0 & 0 & 0 & 0
   & 0 & 0 & 0 & 0 & 0 & 0 & 0 & 0 & 0 & 0 & 0 & 0 & 0 & 0 \\
 0 & 0 & 0 & 0 & 0 & 0 & 0 & 0 & 0 & 0 & 0 & 0 & 0 & 0 & 1 & 1 & 0 & 0 & 0 & 0 & 0 & 0 & 0 & 0 & 0 & 0
   & 0 & 0 & 0 & 0 & 0 & 0 & 0 & 0 & 0 & 0 & 0 & 0 & 0 & 0 \\
\end{array}\right)$}
\newpage

$M_{31}=$\\

\hspace{-4.7cm}
{\tiny $
\left(
\begin{array}{llllllllllllllllllllllllllllllllllllllll}
 0 & 0 & 0 & 0 & 0 & 0 & 1 & 0 & 0 & 0 & 0 & 0 & 0 & 0 & 0 & 0 & 0 & 0 & 0 & 0 & 1 & 0 & 0 & 0 & 0 & 0
   & 0 & 0 & 0 & 0 & 0 & 0 & 0 & 0 & 0 & 0 & 0 & 0 & 0 & 0 \\
 0 & 0 & 0 & 0 & 0 & 0 & 0 & 0 & 0 & 0 & 0 & 0 & 0 & 0 & 0 & 0 & 0 & 0 & 0 & 0 & 0 & 0 & 0 & 0 & 0 & 0
   & 0 & 0 & 0 & 0 & 0 & 0 & 0 & 0 & 0 & 0 & 1 & 1 & 0 & 0 \\
 0 & 0 & 0 & 1 & 0 & 0 & 0 & 0 & 0 & 0 & 0 & 0 & 0 & 0 & 0 & 0 & 0 & 0 & 0 & 0 & 1 & 0 & 0 & 0 & 0 & 0
   & 0 & 0 & 0 & 0 & 0 & 0 & 0 & 0 & 0 & 0 & 0 & 1 & 0 & 0 \\
 0 & 0 & 0 & 0 & 0 & 0 & 0 & 0 & 0 & 0 & 0 & 0 & 0 & 0 & 0 & 0 & 0 & 1 & 0 & 0 & 1 & 0 & 0 & 0 & 0 & 0
   & 0 & 0 & 0 & 0 & 0 & 0 & 0 & 0 & 0 & 0 & 0 & 0 & 0 & 0 \\
 1 & 0 & 0 & 0 & 1 & 0 & 0 & 0 & 0 & 0 & 0 & 0 & 0 & 0 & 0 & 0 & 0 & 0 & 0 & 0 & 0 & 0 & 0 & 0 & 0 & 0
   & 0 & 0 & 0 & 0 & 0 & 0 & 0 & 0 & 0 & 0 & 0 & 0 & 0 & 0 \\
 0 & 0 & 0 & 0 & 0 & 0 & 0 & 0 & 0 & 0 & 0 & 0 & 0 & 0 & 0 & 0 & 0 & 0 & 0 & 0 & 0 & 0 & 0 & 0 & 0 & 0
   & 0 & 0 & 0 & 0 & 0 & 0 & 1 & 0 & 0 & 0 & 0 & 0 & 0 & 0 \\
 0 & 0 & 0 & 1 & 1 & 0 & 0 & 0 & 0 & 0 & 0 & 0 & 0 & 0 & 0 & 0 & 0 & 0 & 0 & 0 & 0 & 0 & 0 & 0 & 0 & 0
   & 0 & 0 & 0 & 0 & 0 & 0 & 0 & 0 & 0 & 0 & 0 & 0 & 0 & 0 \\
 0 & 0 & 0 & 0 & 0 & 0 & 0 & 0 & 0 & 0 & 0 & 0 & 0 & 0 & 0 & 0 & 0 & 0 & 0 & 0 & 0 & 0 & 0 & 0 & 0 & 0
   & 0 & 0 & 0 & 0 & 0 & 0 & 0 & 0 & 0 & 0 & 0 & 1 & 0 & 0 \\
\end{array}
\right)$}
\newpage

$M_{32}=$\\

\hspace{-4.7cm}
{\tiny $ 
\left(
\begin{array}{llllllllllllllllllllllllllllllllllllllll}
 0 & 0 & 0 & 0 & 0 & 0 & 0 & 0 & 0 & 0 & 0 & 0 & 0 & 0 & 0 & 0 & 1 & 1 & 0 & 0 & 0 & 0 & 0 & 0 & 0 & 0
   & 0 & 0 & 0 & 0 & 0 & 0 & 0 & 0 & 0 & 0 & 0 & 0 & 0 & 0 \\
 0 & 0 & 0 & 0 & 0 & 0 & 0 & 0 & 0 & 0 & 0 & 0 & 0 & 0 & 0 & 0 & 1 & 0 & 0 & 0 & 0 & 0 & 0 & 0 & 0 & 0
   & 0 & 0 & 0 & 1 & 0 & 0 & 0 & 0 & 0 & 0 & 0 & 0 & 0 & 0 \\
 0 & 1 & 0 & 0 & 0 & 0 & 0 & 0 & 0 & 0 & 0 & 0 & 0 & 0 & 0 & 0 & 0 & 0 & 0 & 0 & 0 & 0 & 0 & 0 & 0 & 0
   & 0 & 0 & 0 & 0 & 0 & 0 & 0 & 0 & 0 & 0 & 0 & 0 & 0 & 0 \\
 0 & 0 & 0 & 0 & 0 & 0 & 0 & 0 & 0 & 0 & 0 & 0 & 0 & 0 & 0 & 0 & 0 & 0 & 0 & 0 & 0 & 0 & 0 & 0 & 0 & 0
   & 0 & 0 & 0 & 0 & 0 & 0 & 0 & 0 & 0 & 0 & 0 & 1 & 0 & 0 \\
 0 & 0 & 0 & 0 & 0 & 0 & 0 & 0 & 0 & 0 & 0 & 0 & 0 & 0 & 0 & 0 & 0 & 0 & 0 & 0 & 1 & 0 & 0 & 0 & 1 & 0
   & 0 & 0 & 0 & 0 & 0 & 0 & 0 & 0 & 0 & 0 & 0 & 0 & 0 & 0 \\
 0 & 0 & 0 & 0 & 0 & 1 & 0 & 0 & 0 & 0 & 0 & 0 & 0 & 0 & 0 & 0 & 0 & 0 & 0 & 0 & 1 & 1 & 0 & 0 & 0 & 0
   & 0 & 0 & 0 & 0 & 0 & 0 & 0 & 0 & 0 & 0 & 0 & 0 & 0 & 0 \\
 0 & 0 & 0 & 0 & 0 & 0 & 0 & 0 & 0 & 0 & 0 & 0 & 0 & 0 & 0 & 0 & 0 & 0 & 0 & 0 & 0 & 0 & 0 & 0 & 0 & 0
   & 0 & 0 & 0 & 0 & 0 & 0 & 0 & 1 & 0 & 0 & 0 & 0 & 0 & 0 \\
 0 & 0 & 0 & 0 & 0 & 1 & 0 & 0 & 0 & 0 & 0 & 0 & 0 & 0 & 0 & 0 & 0 & 0 & 0 & 0 & 0 & 0 & 0 & 0 & 0 & 0
   & 0 & 0 & 0 & 0 & 0 & 0 & 0 & 1 & 0 & 0 & 0 & 1 & 0 & 0 \\
\end{array}
\right)$}
\newpage

$M_{33}=$\\

\hspace{-4.7cm}
{\tiny $
\left(
\begin{array}{llllllllllllllllllllllllllllllllllllllll}
 0 & 0 & 0 & 0 & 0 & 0 & 0 & 0 & 0 & 0 & 0 & 0 & 0 & 0 & 0 & 0 & 0 & 0 & 1 & 0 & 0 & 0 & 0 & 0 & 0 & 0
   & 0 & 0 & 0 & 0 & 0 & 1 & 0 & 0 & 0 & 0 & 0 & 0 & 0 & 0 \\
 0 & 0 & 0 & 0 & 0 & 0 & 0 & 0 & 0 & 0 & 0 & 0 & 0 & 0 & 0 & 0 & 0 & 0 & 0 & 0 & 0 & 0 & 0 & 0 & 0 & 0
   & 1 & 0 & 0 & 0 & 0 & 1 & 0 & 0 & 0 & 0 & 0 & 0 & 0 & 0 \\
 0 & 0 & 0 & 0 & 0 & 0 & 0 & 0 & 0 & 0 & 0 & 0 & 0 & 0 & 0 & 0 & 0 & 0 & 0 & 0 & 0 & 0 & 0 & 0 & 0 & 0
   & 0 & 0 & 1 & 0 & 0 & 1 & 0 & 0 & 0 & 0 & 0 & 0 & 0 & 0 \\
 0 & 0 & 0 & 0 & 0 & 0 & 0 & 0 & 0 & 0 & 0 & 0 & 0 & 0 & 0 & 1 & 0 & 0 & 0 & 0 & 0 & 0 & 0 & 0 & 0 & 0
   & 0 & 0 & 0 & 0 & 1 & 1 & 0 & 0 & 0 & 0 & 0 & 0 & 0 & 0 \\
 0 & 0 & 0 & 0 & 0 & 0 & 0 & 0 & 0 & 0 & 0 & 0 & 0 & 0 & 0 & 0 & 1 & 0 & 0 & 0 & 0 & 0 & 0 & 0 & 0 & 0
   & 0 & 0 & 0 & 0 & 0 & 0 & 1 & 0 & 0 & 0 & 0 & 0 & 0 & 0 \\
 0 & 0 & 0 & 0 & 0 & 0 & 0 & 0 & 0 & 0 & 0 & 0 & 0 & 0 & 0 & 0 & 0 & 0 & 0 & 0 & 0 & 0 & 0 & 0 & 0 & 1
   & 0 & 0 & 0 & 0 & 0 & 0 & 1 & 0 & 0 & 0 & 0 & 0 & 0 & 0 \\
 0 & 0 & 0 & 0 & 0 & 0 & 0 & 0 & 0 & 0 & 0 & 0 & 0 & 0 & 1 & 0 & 0 & 0 & 0 & 0 & 0 & 0 & 0 & 0 & 0 & 0
   & 0 & 0 & 1 & 0 & 0 & 0 & 1 & 0 & 0 & 0 & 0 & 0 & 0 & 0 \\
 0 & 0 & 0 & 0 & 0 & 0 & 0 & 0 & 0 & 0 & 0 & 0 & 0 & 0 & 0 & 0 & 0 & 0 & 0 & 0 & 0 & 0 & 0 & 0 & 0 & 0
   & 0 & 0 & 0 & 0 & 0 & 1 & 1 & 0 & 0 & 0 & 0 & 0 & 0 & 0 \\
\end{array}
\right)$}

\bibliography{bib}

\end{document}